\documentclass[11pt]{amsart}

\usepackage{amsmath}
\usepackage{amssymb}
\usepackage{amscd}
\usepackage[all]{xy}
\usepackage{sty1}

\makeatletter
\@addtoreset{equation}{section}
\makeatother

\def\bltwo#1{\mathrm{L}^2(#1)}

\def\ccfs#1{\mathrm{A}_{#1}}
\def\fal#1{\mathrm{A}(#1)}
\def\falg{\mathrm{A}(G)}
\def\falh{\mathrm{A}(H)}
\def\faltg{\mathrm{A}_\tau(G)}
\def\falb#1#2{\mathrm{A}_{#1}(#2)} 
\def\falbg#1{\mathrm{A}_{#1}(G)}
\def\falbh#1{\mathrm{A}_{#1}(H)}
\def\falbho#1{\mathrm{A}_{#1}(H_0)}
\def\falby#1{\mathrm{A}_{#1}(Y)}
\def\fsal#1{\mathrm{B}(#1)}
\def\fsalg{\mathrm{B}(G)}
\def\fsalh{\mathrm{B}(H)}
\def\fsalog{\mathrm{B}_0(G)}
\def\fsaltg{\mathrm{B}_\tau(G)}
\def\fsalbg#1{\mathrm{B}_{#1}(G)}
\def\id{\mathrm{id}}

\def\spine#1{\mathrm{A}^*(#1)}
\def\rspine#1{\mathrm{A}^*_0(#1)}
\def\rspineg{\mathrm{A}_0^*(G)}
\def\spineg{\mathrm{A}^*(G)}
\def\rspineh{\mathrm{A}_0^*(H)}
\def\spineh{\mathrm{A}^*(H)}
\def\rspineho{\mathrm{A}_0^*(H_0)}
\def\spineho{\mathrm{A}^*(H_0)}
\def\vng{\mathrm{VN}(G)}
\def\vn#1{\mathrm{VN}_{#1}}
\def\vno#1{\mathrm{VN}(#1)}

\def\wstarg{\mathrm{W}^*(G)}
\def\wcbop#1{\mathcal{CB}^\sigma(#1)}

\def\ctnqg{\wbar{\mathcal{T}}_{nq}(G)}
\def\lutau{\mathcal{LU}_\tau}
\def\rutau{\mathcal{RU}_\tau}
\def\hdg{\mathfrak{HD}(G)}   
\def\rhdg{\mathfrak{HD}_0(G)} 
\def\tg{\mathcal{T}(G)}
\def\thg{\mathcal{T}_0(G)}
\def\tnqg{\mathcal{T}_{nq}(G)}
\def\tnqhg{\mathcal{T}_{nq}^0(G)}
\def\th{\mathcal{T}(H)}
\def\thh{\mathcal{T}_0(H)}
\def\tnqh{\mathcal{T}_{nq}(H)}
\def\tnqhh{\mathcal{T}_{nq}^0(H)}
\def\tho{\mathcal{T}(H_0)}
\def\thho{\mathcal{T}_0(H_0)}
\def\tnqho{\mathcal{T}_{nq}(H_0)}

\def\tnq#1{\mathcal{T}_{nq}(#1)}

\def\mapg{G^{ap}}
\def\tauap{\tau_{ap}}
\def\etaap{\eta_{ap}}
\def\taunq{\tau_{nq}}
\def\mapgtau1{G_{\tau_1}^{ap}}
\def\etatau1ap{\eta^{\tau_1}_{ap}}

\def\cringh{\Omega(H)}
\def\ocringh{\Omega_o(H)}
\def\cbnorm#1{\left\|#1\right\|}

\begin{document}

\newtheorem{qprojlim}{Proposition}[section]
\newtheorem{nonquot}[qprojlim]{Theorem}
\newtheorem{equivtop}[qprojlim]{Lemma}

\newtheorem{spinealgebra}{Proposition}[section]
\newtheorem{subalgebraq}[spinealgebra]{Lemma}
\newtheorem{containment}[spinealgebra]{Proposition}
\newtheorem{intersection}[spinealgebra]{Lemma}
\newtheorem{directsum}[spinealgebra]{Theorem}
\newtheorem{grading}[spinealgebra]{Corollary}

\newtheorem{spectrum}{Theorem}[section]
\newtheorem{spectopology}[spectrum]{Theorem}
\newtheorem{spectopology1}[spectrum]{Corollary}
\newtheorem{hdprincipal}[spectrum]{Proposition}
\newtheorem{rspectrum}[spectrum]{Theorem}
\newtheorem{rspectopology2}[spectrum]{Corollary}
\newtheorem{idempotents}[spectrum]{Proposition}
\newtheorem{whenunital}[spectrum]{Theorem}
\newtheorem{pamaps}[spectrum]{Lemma}
\newtheorem{homomorphism}[spectrum]{Theorem}
\newtheorem{exttop}[spectrum]{Lemma}
\newtheorem{spinerest}[spectrum]{Theorem}
\newtheorem{fsalidem}[spectrum]{Theorem}
\newtheorem{homorange}[spectrum]{Theorem}
\newtheorem{homorange1}[spectrum]{Corollary}

\newtheorem{abdual}{Theorem}[section]
\newtheorem{abdual1}[abdual]{Corollary}

\newtheorem{quotgroup}{Lemma}[section]
\newtheorem{reeandtee}[quotgroup]{Proposition}
\newtheorem{vectorgroup}[quotgroup]{Theorem}
\newtheorem{lattice}[quotgroup]{Corollary}
\newtheorem{mwap}[quotgroup]{Theorem}
\newtheorem{axpb}[quotgroup]{Theorem}
\newtheorem{padic}[quotgroup]{Proposition}
\newtheorem{padicandreal}[quotgroup]{Proposition}

\title[The Spine of a Fourier-Stieltjes Algebra]
{The Spine of a Fourier-Stieltjes Algebra}

\author{Monica Ilie and Nico Spronk}

\maketitle

\footnote{{\it Date}: \today.


We like to thank NSERC for partial funding of our work.}

\begin{abstract}
We define the {\it spine} $\spineg$
of the Fourier-Stieltjes algebra $\fsalg$ of a locally compact
group $G$.  This algebra encodes information about much of the fine
structure of $\fsalg$; particularly information about certain homomorphisms 
and idempotents.

We show that $\spineg$ is graded over a certain semi-lattice,
that of non-quotient locally precompact topologies on $G$.
We compute the spine's spectrum $G^*$, which admits a semi-group
structure.  We discuss homomorphisms from $\spineg$ to $\fsalh$
where $H$ is another locally compact group; and we show that 
$\spineg$ contains the image of every completely bounded homomorphism from 
the Fourier algebra $\falh$ of any amenable group $H$.
We also show that $\spineg$ contains all of the idempotents
in $\fsalg$.
Finally, we compute examples for vector groups, abelian lattices, minimally
almost periodic groups and the $ax+b$-group; and we explore the complexity
of $\spineg$ for the discrete rational numbers and free groups.
\end{abstract}


\section{Introduction}

\subsection{History, Motivation and Plan}

For any locally compact abelian group $G$, Inoue \cite{inoue1} developed
a subalgebra of the measure algebra $\mathrm{M}(G)$.  Let $\tau_G$
be the topology of $G$, and for any locally compact group topology $\tau$
on $G$ which is finer than $\tau_G$, there is a natural isometric homomorphism
from the group algebra $\mathrm{L}^1(G,\tau)$ into $\mathrm{M}(G)$;
let us denote the range $\mathrm{L}^\tau(G)$.
If, for example, $\tau_d$ is the discrete topology, then we obtain the
isomorphism which identifies $\ell^1(G)=\mathrm{L}^1(G,\tau_d)$
with the closed span of the Dirac measures in $\mathrm{M}(G)$.
Then Inoue's algebra is
\begin{equation}\label{eq:inoue}
\mathrm{L}^*(G)=\wbar{\sum_{\tau}\mathrm{L}^\tau(G)}
\end{equation}
where $\tau$ ranges over all locally compact topologies finer than
$\tau_G$ and the overline denotes norm closure.  

In \cite[7.6.5]{taylor}, Taylor identifies for one of his
generalised ``convolution measure algebras" $\mathfrak{M}$,
the closed linear span of all maximal subalgebras which
are isometrically isomorphic to a group algebra of some locally
compact abelian group $G$.  He calls this the {\it spine} of
$\mathfrak{M}$ and denotes it $\mathfrak{M}_0$.  
By \cite[8.2.2]{taylor} (or by
Cohen's homomorphism theorem \cite{cohen} and the
characterisation of the range of homomorphisms from 
abelian group algebras into abelian measure algebras from 
\cite{inoue2}), we see that if $\mathfrak{M}=\mathrm{M}(G)$,
then $\mathfrak{M}_0=\mathrm{L}^*(G)$.

\medskip
In the present article we develop the non-commutative dual
analogue of the spine of an abelian measure algebra.
Thus we let $G$ be a locally compact group and $\fsalg$
denote the {\it Fourier-Stieltjes algebra} of Eymard \cite{eymard}.
(If $G$ is abelian, than $\fsalg\cong\mathrm{M}(\what{G})$ via the
Fourier-Stieltjes transform on the dual group $\what{G}$.)
In this context the {\it Fourier algebra} $\falg$ plays the role
of ``$\mathrm{L}^1(\what{G})$".  Thus the analogue of the spine,
in this context, ought to be the closed span of all of the ``maximal
Fourier subalgebras".

To this end we consider group topologies $\tau$ on $G$ which are coarser
than $\tau_G$, since in the abelian setting the continuous homomorphism
$(G,\tau_G)\to(G,\tau)$ give rise to a dual map $(\hat{G},\hat{\tau})\to
(\hat{G},\tau_{\what{G}})$ so $\hat{\tau}\subseteq\tau_{\what{G}}$.  
We wish to restrict ourselves to such topologies
which will give us a Fourier algebra so we choose
topologies $\tau$ for which the completion of $G$ with respect to $\tau$, 
$G_\tau$, is locally compact.
We call these {\it locally precompact topologies}.  Then there is a continuous
homomorphism with dense range,
$\eta_\tau:G\to G_\tau$, which gives rise to an isometric
homomorphism $u\mapsto u\comp\eta_\tau$ from $\fal{G_\tau}$ to
$\fsalg$.  We denote the range of the homomorphism $\faltg$ and define
the {\it spine} by
\[
\spineg=\wbar{\sum_{\tau\in\tg}\faltg}
\]
where $\tg$ is the family of all locally precompact topologies, coarser
than $\tau_G$.
The notion of maximality of a Fourier subalgebra is suggested to us by
\cite[Theorem 2.1]{bekkakls}; see the proof of Lemma \ref{lem:subalgebraq}.
This motivates us to identify {\it non-quotient topologies} in 
\S\ref{ssec:quotient}.  The existence of such topologies is aided
by appealing to the {\it almost periodic compactification} of $G$.  These
topologies admit a ``supremum" operation, which allows us
to see that $\spineg$ is graded over a semi-lattice.

We note
that the topologies we consider are not necessarily Hausdorff on $G$,
in particular when the group fails to be maximally almost periodic.  We
will deal specially with the Hausdorff topologies in the {\it reduced spine},
$\rspineg$.  We show in \S \ref{ssec:rspinespec}
that $\rspineg$ is unital, if and only if $G$ is 
maximally almost periodic, if and only if $\rspineg=\spineg$.

Motivated by the papers of Inoue \cite{inoue1,inoue2}, 
in Section \ref{sec:charhomo} we 
exploit the semi-lattice structure on our locally precompact non-quotient
topologies to construct
the spectrum $G^*$ of $\spineg$.  Using $G^*$ we can extend
our generalisation \cite{ilies} of Cohen's homomorphism theorem
\cite{cohen} to study homomorphisms from $\spineg$ to $\fsalh$
for another locally compact group $H$.  
We are forced, however, to consider the natural operator space structures
on $\spineg$ and $\fsalh$, and as well restrict ourselves to amenable $G$.
Thus we obtain a characterisation of
{\it completely bounded} homomorphisms from the spine $\spineg$ for an
amenable group $G$, to $\fsalh$ for a general locally compact group $H$.
Moreover, we show that $\spineh$ is exactly the subspace of $\fsalh$ which
contains the images of all such homomorphisms.  We also obtain the facts
that if $H_0$ is an open subgroup of $H$, then $\spineh|_{H_0}=
\spineho$, and that $\spineh$ contains all of the idempotents from
$\fsalh$.

In Section \ref{sec:abgroup} we show that for an abelian locally 
compact group $G$,
that our non-quotient locally precompact group topologies which 
are coarser than
$\tau_G$, are in a natural duality with those locally compact group 
topologies on
$\what{G}$ which are finer that $\tau_{\what{G}}$.  This shows that
$\spineg\cong\mathrm{L}^*(\what{G})$, and that our results generalise those
of Inoue.

In the last section we compute some examples.  For some Lie groups 
we compute examples explicitly.  We even do explicit computations
for some abelian groups, including vector groups, which we have not
found in the existing literature.  Here we use the fact that $G^*$ is a 
semi-topological semigroup, which we dub the {\it spine compactification}
of $G$, to obtain some interesting semi-groups.  We merely scratch the surface
of the problem of attempting to compute the spine for certain algebraic groups
and for free groups.

As this article was being written, some applications of the ideas
discovered have started new projects.  First of them is that 
the present authors can determine the structure
of the subalgebra of $\fsalg$ which is generated by the idempotents.
Second of them is by V.\ Runde and the second named author:  we 
have characterised some amenability properties of $\spineg$.  Moreover, in this
investigation, some examples
of non-compact groups $G$ have been discovered for which $\fsalg$
is both operator amenable and weakly amenable, contrary to standing
conjectures.  Both shall appear in forthcoming articles.

\subsection{Notation}

If $G$ is any locally compact group let $\fsalg$ denote its {\it 
Fourier-Stieltjes algebra} and $\falg$ denote its {\it Fourier algebra}, as
defined in \cite{eymard}.  We recall that $\fsalg$ consists of all
matrix coefficients of continuous unitary representations, i.e.\
functions of the form $s\mapsto\inprod{\pi(s)\xi}{\eta}$ where
$\pi:G\to\fU(\fH)$ is a homomorphism, continuous when the unitary group
$\fU(\fH)$ on the Hilbert space $\fH$ is endowed with the weak operator 
topology.  We also recall that $\falg$ is the space of all matrix
coefficients of the left regular representation $\lam_G:G\to\fU(\bltwo{G})$,
given by left translation operators on
$\bltwo{G}$, the Hilbert space of (equivalence classes of)
square-integrable functions. 
The Fourier-Stieltjes algebra
is the predual of the {\it enveloping von Neumann algebra} $\wstarg$ which
is generated by the universal representation $\varpi_G$ \cite{dixmierB}.
This norm makes $\fsalg$ into a Banach algebra.
If $\pi:G\to\fU(\fH)$ is any continuous unitary representation of $G$, we let
$\vn{\pi}$ be the von Neumann algebra generated by $\pi(G)$ and
$\ccfs{\pi}$ be the closed linear subspace generated by matrix coefficients
$\{\inprod{\pi(\cdot)\xi}{\eta}:\xi,\eta\in\fH\}$.  We have that
the annihilator $\ccfs{\pi}^\perp$ in $\wstarg$ is closed ideal
for which $\wstarg/\ccfs{\pi}^\perp\cong\vn{\pi}$ as von Neumann algebras,
whence we have isometric duality ${\ccfs{\pi}}^*\cong\vn{\pi}$.  
In particular, $\falg$ is the predual
of the {\it group von Neumann algebra} $\vng=\vn{\lam_G}$.  We note that
$\falg$ is an ideal in $\fsalg$ and is semi-simple with character space
bicontinuously isomorphic to $G$.

We will make enough use of operator spaces that a few words are in order.
Our standard reference on this topic is \cite{effrosrB}, though we will 
indicate other references below.  
Any C*-algebra $\fA$ is an {\it operator space}
in the sense that for $n=1,2,\dots$ the algebra of $n\cross n$ matrices
over $\fA$, $\mathrm{M}_n(\fA)$ admits a unique norm which makes it into a 
C*-algebra.  A linear map $T:\fA\to\fB$ between C*-algebras is called
{\it completely bounded} if it is bounded and 
its amplifications $T^{(n)}:\mathrm{M}_n(\fA)\to
\mathrm{M}_n(\fB)$, given by $T^{(n)}[a_{ij}]=[Ta_{ij}]$, give a bounded
family of norms $\left\{\norm{T^{(n)}}:n=1,2,\dots\right\}$.  In this case
we write $\cbnorm{T}=\sup\left\{\norm{T^{(n)}}:n=1,2,\dots\right\}$.

If $\fM$ and $\fN$ are von Neumann algebras with
preduals $\fM_*$ and $\fN_*$, we say a map $\Phi:\fM_*\to\fN_*$ is
{\it completely bounded}, if its adjoint 
$\Phi^*:\fN\to\fM$ is such.  However, it is often convenient to
consider completely bounded maps on the space $\fM_*$ by noting that
it admits an operator space structure via the identifications
$\mathrm{M}_n(\fM_*)\cong\wcbop{\fM,\mathrm{M}_n}$, $n=1,2,\dots$, where 
$\wcbop{\fM,\mathrm{M}_n}$ is the space of normal completely bounded maps 
from $\fM$ to the finite dimensional von Neumann algebra of $n\cross n$
complex matrices \cite{blecher,effrosr}.  We note that these
spaces $\fM_*$, with the above matricial structures,
are completely isometrically isomorphic to subspaces of C*-algebras 
\cite{ruan0}, which are not generally operator subalgebras.
In particular, these
structures are used to create the {\it operator projective tensor product}
$\fM_*\what{\otimes}\fN_*$ \cite{blecherp,effrosr}.  This tensor product
admits the very useful formula 
$(\fM_*\what{\otimes}\fN_*)^*\cong\fM\wbar{\otimes}\fN$
\cite{effrosr}, where $\fM\wbar{\otimes}\fN$ is the von Neumann tensor product.

\section{The Lattice of Locally Precompact Group Topologies}

\subsection{Locally Precompact Topologies}

Let $G$ be a locally compact group.  We will denote by $\tau_G$
the topology on $G$.  As is conventional in the literature, we will
suppose that $\tau_G$ is Hausdorff.  If $H$ is another locally compact group
with topology $\tau_H$, and $\eta:G\to H$ is a continuous homomorphism,
then
\begin{equation}\label{eq:indtop}
\eta^{-1}(\tau_H)=\{\eta^{-1}(U):U\in\tau_H\}
\end{equation}
forms a group topology on $G$.  We note that $\eta^{-1}(\tau_H)$ is coarser
than $\tau_G$ and is Hausdorff exactly when $\eta$ is injective.  Letting
$\fP(G)$ denote the power set of $G$ we let
\begin{equation}\label{eq:deftg}
\tg=\left\{\tau\subseteq\fP(G):
\begin{matrix}\text{there exist a locally compact group }H\text{ and} \\
\text{a continuous homomorphism }\eta:G\to H \\ \text{such that }
\tau=\eta^{-1}(\tau_H) \end{matrix}\right\}.
\end{equation}
We call the elements of $\tg$ {\it locally precompact} topologies on $G$.

If $\tau\in\tg$ and $\tau=\eta^{-1}(\tau_H)$, as in (\ref{eq:indtop}),
we call the pair $(\eta,H)$ a {\it representation} of $\tau$.
If $(\eta_j,H_j)$, $j=1,2$, are each representations of $\tau$, then
there is a topological isomorphism 
\[
\iota_2^1:\wbar{\eta_1(G)}\to\wbar{\eta_2(G)}\text{ such that }
\iota_2^1\comp\eta_1=\eta_2.
\]
Indeed, if $\fU_\tau$ is a neighbourhood basis for $\tau$ at the identity
$e$ in $G$, then we see that
\[
N_\tau=\bigcap_{U\in\fU_\tau}U
\]
coincides with each $\ker\eta_j$, so $N_\tau$ is a $\tau$-closed
normal subgroup of $G$.  Let $q_\tau:G\to G/N_\tau$ be the quotient map.
The first isomorphism theorem provides us
with isomorphisms $\iota_j:G/N_\tau\to\eta_j(G)$, $j=1,2$, for which
$\iota_j\comp q_\tau=\eta_j$ and 
the topologies $\iota_j^{-1}(\tau_{H_j})$ each coincide with
the quotient topology $q_\tau(\tau)=\{q_\tau(U):U\in\tau\}$.
Thus we obtain an topological isomorphism $\iota_2^1:\eta_1(G)\to\eta_2(G)$,
which by uniform continuity (see \cite[5.40(a)]{hewittrI}, 
and \cite[6.26]{kelley}), 
extends uniquely to a continuous isomorphism
$\bar{\iota}_2^1:\wbar{\eta_1(G)}\to\wbar{\eta_2(G)}$.  Changing the order
we obtain a continuous isomorphism 
$\bar{\iota}_1^2:\wbar{\eta_2(G)}\to\wbar{\eta_1(G)}$ which is the inverse
of $\bar{\iota}_2^1$, hence $\bar{\iota}_2^1$ is a topological isomorphism,
which we will denote $\iota_2^1$.  Hence, if $\tau\in\tg$ we may select
a representation $(\eta_\tau,G_\tau)$, which is unique up to isomorphism,
and satisfies

{\bf (i)} $\eta_\tau:G\to G_\tau$ is a continuous homomorphism with
dense range,

{\bf (ii)} $\ker\eta_\tau=N_\tau$, and

{\bf (iii)} $\tau=\eta_\tau^{-1}(\tau_{G_\tau})$.

\noindent We will call $G_\tau$ the {\it completion} of $G$ with 
respect to $\tau$,
and $N_\tau$ the {\it kernel} of $\tau$.  We will let $q_\tau:
G\to G/N_\tau$ be the quotient map, and $\iota_\tau:G/N_\tau\to G_\tau$
be the continuous homomorphism
provided by the first isomorphism theorem.  We note that
$\iota_\tau$ is always injective, though it may not be surjective.

If $\tau_1\subseteq\tau_2$ in $\tg$, then there is a unique homomorphism
\[
\eta_{\tau_1}^{\tau_2}:G_{\tau_2}\to G_{\tau_1}
\text{ such that }\eta_{\tau_1}^{\tau_2}\comp\eta_{\tau_2}=\eta_{\tau_1}.
\]
Indeed, we have that
\[
N_{\tau_1}=\bigcap_{U\in\fU_{\tau_1}}U\supseteq
\bigcap_{V\in\fU_{\tau_2}}V=N_{\tau_2}
\]
and hence the second isomorphism theorem provides us with a continuous
homomorphism $q_{\tau_1}^{\tau_2}:G/N_{\tau_1}\to G/N_{\tau_2}$
such that $q_{\tau_1}^{\tau_2}\comp q_{\tau_2}=q_{\tau_1}$.
This induces a continuous homomorphism $\til{q}_{\tau_1}^{\tau_2}:
\iota_{\tau_1}(G/N_{\tau_1})\to\iota_{\tau_2}(G/N_{\tau_2})\subseteq
G_{\tau_2}$, which extends uniquely to a continuous homomorphism
$\eta_{\tau_1}^{\tau_2}:G_{\tau_2}\to G_{\tau_1}$, as required.
We note that if we have a third topology $\tau_0$ in $\tg$, 
$\tau_0\subseteq\tau_1$, then we obtain the relation
\begin{equation}\label{eq:etarel}
\eta_{\tau_0}^{\tau_1}\comp\eta_{\tau_1}^{\tau_2}=
\eta_{\tau_0}^{\tau_2}.
\end{equation}
In the case that $\tau_2=\tau_G$ we recover the relation
$\eta_{\tau_0}^{\tau_1}\comp\eta_{\tau_1}=\eta_{\tau_0}$.

It will be convenient for us to distinguish
\begin{equation}\label{eq:defthg}
\thg=\{\tau\in\tg:\tau\text{ is Hausdorff}\}
\end{equation}
which is exactly the set of those topologies $\tau$ in $\tg$ for which
$N_\tau=\{e\}$.

\bigskip
We note that we can develop $\thg$ by an alternate method.
If $\tau$ is any Hausdorff topology on $G$ which is coarser than $\tau_G$
we define the {\it left uniformity} on $G$ with respect to $\tau$ by
\[
\lutau=\{W\subseteq G\cross G:W\supseteq m_l^{-1}(U)
\text{ for some }U\in\fU_\tau\}
\]
where $m_l:G\cross G\to G$ is given by $m_l(s,t)=s^{-1}t$.
We may similarly define the right uniformity on $G$ with respect to $\tau$,
$\rutau$.  We say that
$\tau$ is {\it locally precompact} if there exists a completion
of $G$ with respect to $\lutau$, $G_\tau$, which is a locally compact group.  
By \cite[9.26 and 10.15]{roelcked},
the existence of such a completion with respect to $\lutau$
is equivalent to having a locally compact completion with respect to $\rutau$,
and these two completions are topologically isomorphic.
Note that if $\tau$ is a locally compact topology, then $G$ is already
complete with respect to  $\lutau$, and we have $G=G_\tau$.

\subsection{The Lattice Operation}\label{subsec:lattice}
If $\tau_1,\tau_2\in\tg$, we let $\tau_1\vee\tau_2$ denote the coarsest
topology which makes each map $\eta_{\tau_j}:G\to G_{\tau_j}$, $j=1,2$
continuous.  It is clear that $\tau_1\vee\tau_2$ has
base $\tau_1\cup\tau_2$, and is thus coarser than $\tau_G$.  
Moreover, $\tau_1\vee\tau_2$ is locally precompact for if $\del:G\to
G_{\tau_1}\cross G_{\tau_2}$ is the map given by 
$\del(s)=(\eta_{\tau_1}(s),\eta_{\tau_2}(s))$, then 
\[
\tau_1\vee\tau_2=\del^{-1}(\tau_1\cross\tau_2)
=\{\del^{-1}(W):W\in\tau_1\cross\tau_2\}
\]
where $\tau_1\cross\tau_2$ is the product topology on 
$G_{\tau_1}\cross G_{\tau_2}$.  Thus we have that
\begin{equation}\label{eq:supmap}
\eta_{\tau_1\vee\tau_2}:G\to G_{\tau_1\vee\tau_2}\subset 
G_{\tau_1}\cross G_{\tau_2}:s\mapsto(\eta_{\tau_1}(s),\eta_{\tau_2}(s)).
\end{equation}

We note some very basic relations.  Let $\tau_1,\tau_2,\tau_3\in\tg$.  Then
\begin{equation}\label{eq:taurel}
\tau_1\vee\tau_2=\tau_2\vee\tau_1\quad\aand\quad \tau_1\vee\tau_2=\tau_2
\iff \tau_1\subseteq\tau_2.
\end{equation}
Moreover, the lattice operation is associative:
$(\tau_1\vee\tau_2)\vee\tau_3=\tau_1\vee(\tau_2\vee\tau_3)$.
Also, since $N_{\tau_1\vee\tau_2}=N_{\tau_1}\cap N_{\tau_2}$, it follows
that
\begin{equation}\label{eq:thgideal}
\text{if }\tau\in\tg\text{ and }\tau_0\in\thg,\text{ then }
\tau\vee\tau_0\in\thg\text{ too.}
\end{equation}

In order to describe the lattice operation on infinite families,
let us first
note that $\tg$ is a directed set via inclusion: any pair
$\tau_1,\tau_2$ of elements is dominated by $\tau_1\vee\tau_2$.  
Any directed subset
$\fS$ of $\tg$ thus gives rise to an inverse mapping system:
\[
\left\{G_\tau,\eta^{\tau_2}_{\tau_1}:\tau\in\fS,
\tau_1\subseteq\tau_2\iin\fS\right\}
\]
by the relations (\ref{eq:etarel}).  Thus we obtain the projective limit
\begin{equation}\label{eq:projlim}
G_\fS=\underset{\tau\in\fS}{\underleftarrow{\lim}}G_\tau
=\left\{(s_\tau)\in\prod_{\tau\in\fS}G_\tau:\eta^{\tau_2}_{\tau_1}(s_{\tau_2})
=s_{\tau_1}\iif \tau_1\subseteq\tau_2\iin\fS\right\}.
\end{equation}
We denote the topology on this group by
$\tau_\fS=\bigvee_{\tau\in \fS}\tau$.
It is well-known that $G_\fS$ is a complete topological group (see
\cite[Section 25]{husain}, for example).  However, it is not clear that
$G_\fS$ is locally compact, in general.  The following special
case seems known (there is such a comment on page 133 of
\cite{greene}, for example).  Since
we could not find a proof in any of the standard references \cite{hewittrI},
\cite{husain} or \cite{montgomeryz}, we will supply a brief one.

We say that a continuous map $\vphi:X\to Y$, where $X$, $Y$ are Hausdorff
spaces, is {\it proper}, if $\vphi^{-1}(K)$ is compact in $X$
for every compact subset $K$ of $Y$.

\begin{qprojlim}\label{prop:qprojlim}
Let $\fQ$ be a directed subset of $\tg$ which enjoys the property that if
$\tau_1\subseteq\tau_2$ in $\fQ$, then $\eta^{\tau_2}_{\tau_1}:G_{\tau_1}
\to G_{\tau_2}$ is a proper map.
Then the projective limit $G_\fQ$ is
locally compact.  
\end{qprojlim}

\proof  A subbase for $\tau_\fQ$ is formed by the sets
\begin{equation*}
\left\{\left(\prod_{\tau\in\fQ\setdif\{\tau_0\}} G_\tau\times U_{\tau_0}\right)
\cap G_\fQ
:\begin{matrix} \tau_0\in\fQ\aand U_{\tau_0}\subset G_{\tau_0}
\text{ is a} \\ \text{relatively compact open set}\end{matrix}\right\}.
\end{equation*}
Let us fix one of these neighbourhoods.  If $\tau\in\fQ\setdif\{\tau_0\}$ 
then either $\tau\supset\tau_0$ or there exists a $\tau'$ in $\fQ$
such that $\tau'\supset\tau_0$ and $\tau'\supset\tau$.  Thus we may
express the closure
\begin{align*}
\Biggl(\prod_{\tau\in\fQ\setdif\{\tau_0\}}  & G_\tau
\times \wbar{U}_{\tau_0}\Biggr)
\cap G_\fQ   \\
=&\left(\prod_{\tau\supset\tau_0}(\eta^\tau_{\tau_0})^{-1}(\wbar{U}_{\tau_0})
\times\wbar{U}_{\tau_0}
\times\prod_{\tau\not\supseteq\tau_0}
\eta^{\tau'}_\tau
\bigl((\eta^{\tau'}_{\tau_0})^{-1}(\wbar{U}_{\tau_0})\bigr)
\right)\cap G_\fQ. 
\end{align*}
The latter set is clearly compact. \endpf

If $\{\tau_j\}_{j\in J}\subset\tg$, then we let $\fF$ be the collection
of finite subsets of $J$ and let for $F\iin\fF$, $\tau_F=
\bigvee_{j\in F}\tau_j$ which is in $\tg$.  Then we may realise
$\bigvee_{j\in J}\tau_j=\bigvee_{F\in\fF}\tau_F$ as a projective
limit topology.  Thus it follows Proposition \ref{prop:qprojlim}
that $\bigvee_{j\in J}\tau_j$ is locally precompact if
for each pair of indices $j_1,j_2$ we have that 
$\eta^{\tau_{j_1}\vee\tau_{j_2}}_{\tau_{j_1}}:
G_{\tau_{j_1}\vee\tau_{j_2}}\to G_{\tau_{j_1}}$ is a proper map.

See \S \ref{subsec:dual} for a notion of $\tau_1\wedge\tau_2$.

\subsection{Quotient Topologies}\label{ssec:quotient}
If $\tau,\tau_1\in\tg$ we sat that $\tau$ is a {\it quotient} of $\tau_1$ if

{\bf (i)}  $\tau\subset\tau_1$, and

{\bf (ii)} $\eta^{\tau_1}_\tau:G_{\tau_1}\to G_\tau$ is a 
proper map.

\noindent Recall that proper maps are defined just before Proposition
\ref{prop:qprojlim}.  We note that condition (ii) is equivalent to 

{\bf (ii')} $\ker\eta^{\tau_1}_\tau$ is compact and $G_\tau\cong
G_{\tau_1}/\ker\eta^{\tau_1}_\tau$, homeomorphically.

\noindent Indeed, proper maps are closed (see \cite[Corollary 3.11]{ilies})
so $\eta^{\tau_1}_\tau(G_{\tau_1})=G_\tau$.  Also
$\ker\eta^{\tau_1}_\tau=(\eta^{\tau_1}_\tau)^{-1}(\{e\})$ is compact.
Finally, $\eta^{\tau_1}_\tau$ is open for if $U$ is a relatively compact
open subset of $G_{\tau_1}$, then $V=\eta^{\tau_1}_\tau(U)$ is relatively 
compact in $G_\tau$ so there is a compact neighbourhood $C$ of $V$.
We find that $C\setdif V$ is closed, whence $V$ is open.  
Indeed, if $(s_i)$ is a net in
$C\setdif V$ converging to $s_0$, we let $(t_i)$ be a net from
$(\eta^{\tau_1}_\tau)^{-1}(C)$ for which $\eta^{\tau_1}_\tau(t_i)=s_i$. Then if
$t_0$ is any cluster point of $(t_i)$, we have
that $\eta^{\tau_1}_\tau(t_0)=s_0$
and $t_0\not\in U$, so $s_0\not\in V$, i.e.\ $s_0\in C\setdif V$.

Let
\begin{equation}\label{eq:tnqg}
\tnqg=\left\{\tau\in\tg:\begin{matrix}
\text{ there is no }\tau_1\in\tg\text{ such that} \\
\tau\text{ is a non-trivial quotient of }\tau_1\end{matrix}\right\}.
\end{equation}
We call the elements of $\tnqg$ {\it non-quotient} topologies.
We note that $\tau_G\in\tnqg$.
In order the characterise non-quotient topologies amongst
locally precompact topologies we will make use of a special topology
defined below.

A {\it compactification} of $G$ is a pair $(\eta,S)$ where
$S$ is a semi-topological semigroup and $\eta:G\to S$ is
a continuous homomorphism with dense range.  See the book
\cite{berglundjm} for more on this.  If $\tau\in\tg$ and
$G_\tau$ is compact, then $(\eta_\tau,G_\tau)$ is a group
compactification of $G$, and we say that $\tau$ is {\it precompact}.
The {\it almost periodic compactification} of $G$,
$(\etaap,\mapg)$ is maximal amongst all group compactifications
in the sense that every locally precompact topology $\tau$ on $G$
is a quotient of $\tauap=\etaap^{-1}(\tau_{\mapg})$.

\begin{nonquot}\label{theo:nonquot}
Let $\tau\in\tg$.  Then the following all hold.

{\bf (i)}  The topology $\taunq=\tau\vee\tauap$ admits $\tau$ as
a quotient topology.

{\bf (ii)}  If $\tau$ is a quotient of another topology $\tau_1$ in $\tg$, then
$\tau_1$, in turn, is a quotient of $\taunq$.  Hence $\taunq$ is the unique
non-quotient which admits $\tau$ as a quotient.

{\bf (iii)} $\tau\in\tnqg\quad\iff\quad\tau=\tau\vee\tauap
\quad\iff\quad\tau\supseteq\tauap$.
\end{nonquot}

Let us begin with a lemma.

\begin{equivtop}\label{lem:equivtop}
Let $H$ be a locally compact group, $K$ be a compact normal subgroup
with quotient map $q_K:H\to H/K$,
and $(\eta,H_1)$ be a group compactification of $H$ with
$K\cap\ker\eta=\{e_H\}$.  Then the homomorphism
$s\mapsto\bigl(q_K(s),\eta(s)\bigr)$ from $H$ to $H/K\,\times H_1$
is injective and bicontinuous, thus a homeomorphism onto its range. 
\end{equivtop}

\proof  If for $s,s_1\iin H$, $\bigl(q_K(s),\eta(s)\bigr)=
\bigl(q_K(s_1),\eta(s_1)\bigr)$ then $s_1=sk$ for some $k\iin K$
and hence $\eta(s_1)=\eta(s)\eta(k)$, so $k\in K\cap\ker\eta=\{e_H\}$.
Thus the map is injective.  It is clearly continuous.

Let $(s_i)$ be a net in $H$ and $s_0\iin H$ be such that
$\lim_i\bigl(q_K(s_i),\eta(s_i)\bigr)=\bigl(q_K(s_0),\eta(s_0)\bigr)$.
Then for any relatively compact neighbourhood $U$ of $s_0$ we have
that $(s_i)$ is eventually in $UK$.  We then have that $(s_i)$
has a cluster point $s_1$ in each such $UK$.  Since $U$ is arbitrary,
$s_1$ must be in $s_0K$, whence $s_1=s_0k$ for some $k\iin K$.
If $(s_{i'})$ is a subnet converging to $s_1$ then we see that
\[
\eta(s_0)\eta(k)=\eta(s_1)=\lim_{i'}\eta(s_{i'})=\lim_i\eta(s_i)
=\eta(s_0)
\]
so $k\in K\cap\ker\eta=\{e_H\}$.  Thus $s_1=s_0$ is the only cluster point
of $(s_i)$ and we see that $\lim_i s_i=s_0$.  Thus
$\bigl(q_K(s),\eta(s)\bigr)\mapsto s$ is continuous. \endpf

\noindent {\bf Proof of Theorem \ref{theo:nonquot}.}
Let $\tau\in\tg$.  First, let us see that $\tau$ is a quotient of
$\tau\vee\tauap$.  Let $\eps=\{\varnothing,G\}$ be
the trivial topology, which is precompact as it is
induced by the trivial homomorphism $\eta_\eps:G\to 1$.
Then $\tau=\tau\vee\eps$ and $\tau\vee\eps$ is 
a quotient of $\tau\vee\tauap$ since
$\eta^{\tau\vee\tauap}_{\tau\vee\eps}$ is the restriction,
to $G_{\tau\vee\tauap}$, 
of the quotient map $\id_{G_\tau}\cross\eta^{\tauap}_\eps:
G_\tau\cross\mapg\to G_\tau\cross 1$.

Let $\tau_1$ be any topology of which $\tau$ is a quotient.
Let $(\etatau1ap,\mapgtau1)$ denote the almost periodic compactification
of $G_{\tau_1}$, and $\tau_{1,ap}$ the precompact topology on $G$
induced by $\etatau1ap\comp\eta_{\tau_1}$.  Then by Lemma
\ref{lem:equivtop}, above, 
$s\mapsto\bigl(\eta^{\tau_1}_\tau(s),\etatau1ap(s)\bigr):
G_{\tau_1}\to G_\tau\cross\mapgtau1$ is a bicontinuous isomorphism,
so $\tau_1=\tau\vee\tau_{1,ap}$.  However, $\tau\vee\tau_{1,ap}$
is clearly a quotient of $\taunq=\tau\vee\tauap$.  
Thus (ii) holds.  Moreover, if
$\tau_2$ is any topology of which $\taunq$ is a quotient, then
$\tau$ is also a quotient of $\tau_2$, and hence $\tau_2$ is a quotient
of $\taunq$, so $\tau_2=\taunq$.  Thus $\taunq\in\tnqg$ and
is the unique such topology which has $\tau$ as a quotient; whence we
obtain (i).

Part (iii) is now just a straightforward
application of the lattice relations (\ref{eq:taurel}). \endpf

\subsection{Semilattice Structure on $\tnqg$}\label{subsec:latsemigrp}

A {\it semilattice} consists of a set $\fS$ and
an associative binary operation $\fS\cross \fS\to \fS:
(\sig_1,\sig_2)\mapsto \sig_1\mult \sig_2$
for which $\sig_1\mult \sig_2=\sig_2\mult\sig_1$ and $\sig\mult\sig=\sig$ 
for every $\sig_1,\sig_2,\sig\iin\fS$.  Thus a semilattice is
a commutative idempotent semigroup.

We note that $\tg$ is a semilattice under the operation $(\tau_1,
\tau_2)\mapsto\tau_1\vee\tau_2$ by (\ref{eq:taurel}).
We note that $\tnqg$ is a sub-semilattice by Theorem \ref{theo:nonquot} (iii),
since for $\tau_1,\tau_2\iin\tnqg$ we have that
\[
(\tau_1\vee\tau_2)\vee\tauap=(\tau_1\vee\tauap)\vee(\tau_2\vee\tauap)
=\tau_1\vee\tau_2.
\]

The topology $\tau_G\in\tnqg$, so (\ref{eq:taurel}) tells us that
\begin{equation}\label{eq:taugmax}
\tau_G\vee\tau=\tau_G\text{ for any }\tau\iin\tnqg.  
\end{equation}
Thus $\{\tau_G\}$ is an {\it ideal} in $\tnqg$.
Moreover, the almost periodic topology $\tauap$ is the unit
for $\tnqg$ by Theorem \ref{theo:nonquot} (iii).

We note that $\tnqhg$ is an ideal in
$\tnqg$ by (\ref{eq:thgideal}).
Thus $G$ is a {\it maximally almost periodic} 
group, i.e.\ $\etaap:G\to\mapg$ is injective so $\tauap\in\thg$, exactly when
$\tnqg=\tnqhg$.  This fact will be used in Corollary \ref{cor:grading} and
Theorem \ref{theo:whenunital}.

\section{The Spine $\spineg$ }

\subsection{Definition of the Spine}
If $\tau\in\tg$ we let
\begin{equation}\label{eq:faltg}
\faltg=\fal{G_\tau}\comp\eta_\tau\subseteq\fsalg.
\end{equation}
Since $\eta_\tau:G\to G_\tau$ is a homomorphism with dense range,
$\faltg$ is a closed subalgebra of $\fsalg$ which is
isometrically isomorphic to $\fal{G_\tau}$ by \cite[2.10]{arsac}.
To be more precise, if $\lam_\tau=\lam_{G_\tau}\comp\eta_\tau:
G\to\fU(\bltwo{G_\tau})$, where $\lam_{G_\tau}$ is the left regular 
representation of $G_\tau$, then
\[
\faltg=\ccfs{\lam_\tau}
\]
the space of coefficient functions for $\lam_\tau$.
Thus $u\mapsto u\comp\eta_\tau:\fal{G_\tau}\to\ccfs{\lam_\tau}$ is a complete 
isometry, for its adjoint map $\vn{\lam_\tau}\to\vno{G_\tau}$ 
is a $*$-isomorphism.

We let
\begin{equation}\label{eq:spinedef}
\spineg=\wbar{\sum_{\tau\in\tg}\faltg}.
\end{equation}
and call this subspace the {\it spine} of $\fsalg$.  Moreover, we let
\begin{equation}\label{eq:rspinedef} 
\rspineg=\wbar{\sum_{\tau\in\thg}\faltg}.
\end{equation}
and call it the {\it reduced spine}.

\subsection{Multiplicative Properties}
Let us first establish that $\spineg$ is an algebra.  This is
an extension of the fact that each $\faltg\cong\fal{G_\tau}$ is
itself an algebra.

\begin{spinealgebra}\label{prop:spinealgebra}
If $\tau_1,\tau_2\in\tg$ then $\falbg{\tau_1}\falbg{\tau_2}\subseteq
\falbg{\tau_1\vee\tau_2}$.  In fact
\[
\wbar{\spn}\falbg{\tau_1}\falbg{\tau_2}=\falbg{\tau_1\vee\tau_2}.
\]
\end{spinealgebra}

\proof  Let $u_j\in\falbg{\tau_j}$ for $j=1,2$.
Since $\falbg{\tau_j}=\ccfs{\lam_{\tau_j}}$ 
there are $\xi_j,\eta_j\iin\bltwo{G_{\tau_j}}$ such that
$u_j=\inprod{\lam_{\tau_j}(\cdot)\xi_j}{\eta_j}$.  Then
\[
u_1u_2=\inprod{\lam_{\tau_1}\otimes\lam_{\tau_2}(\cdot)
\xi_1\otimes\xi_2}{\eta_2\otimes\eta_2}.
\]
Since $\bltwo{G_{\tau_1}}\otimes^2\bltwo{G_{\tau_2}}
\cong\bltwo{G_{\tau_1}\cross G_{\tau_2}}$, we have a unitary equivalence
$\lam_{\tau_1}\otimes\lam_{\tau_2}\cong\lam_{\tau_1\vee\tau_2}$.
Thus $\falbg{\tau_1}\falbg{\tau_2}\subseteq\falbg{\tau_1\vee\tau_2}$.
By the operator projective tensor product formula of \cite{effrosr}
we obtain (completely) isometric identifications
\[
\falbg{\tau_1}\what{\otimes}\falbg{\tau_2}
\cong\fal{G_{\tau_1}}\what{\otimes}\fal{G_{\tau_2}}
\cong\fal{G_{\tau_1}\cross G_{\tau_2}}.
\]
The multiplication map $\falbg{\tau_1}\what{\otimes}\falbg{\tau_2}
\to\falbg{\tau_1\vee\tau_2}$ is thus isomorphic to the restriction
map $w\mapsto w|_{G_{\tau_1\vee\tau_2}}:\fal{G_{\tau_1}\cross G_{\tau_2}}
\to\fal{G_{\tau_1\vee\tau_2}}$.  However, the restriction map is
surjective by \cite[Theorem 3]{takesakit}, or \cite{herz}. \endpf

Recall that $\fsalg$ has an involution on it given by pointwise complex 
conjugation.  Thus a $*$-subspace of $\fsalg$ is any subspace closed
under this involution.  

\begin{subalgebraq}\label{lem:subalgebraq}
Let $\tau\in\tg$.  A subset $\fA$ of $\faltg$ is a 
translation invariant, closed $*$-subalgebra 
if, and only if, there exists a
quotient topology $\tau_0$ of $\tau$ in $\tg$ such that $\fA=\falbg{\tau_0}$.
Moreover, if $\fA$ is point separating on $G$, i.e.\ $\tau\in\thg$,
then $\tau_0\in\thg$.
\end{subalgebraq}

\proof If $\fA$ is  translation invariant, closed $*$-subalgebra of $\faltg$, 
then there must be a corresponding translation invariant, closed 
$*$-subalgebra of $\fal{G_\tau}$.
By \cite[Theorem 2.1]{bekkakls}, every translation invariant,
closed $*$-subalgebra of $\fal{G_\tau}$ is of
the form $\fal{G_\tau\!\!:\!\! K}$, the subalgebra of functions which are
constant on cosets of a compact normal subgroup $K$.  
The composition of $\eta_\tau$
with the quotient map $G_\tau\to G_\tau/K$ gives a topology $\tau_0$
which is a quotient of $\tau$.  Moreover, we obtain the identifications 
\[
\fA=\fal{G_\tau\!:\! K}\comp\eta_\tau=
\fal{G_\tau/K}\comp\eta_{\tau_0}=\falbg{\tau_0}.
\]
If $\fA$ is point separating, then we see that $\eta_\tau(G)\cap K=\{e\}$
so $N_{\tau_0}=\{e\}$, so $\tau_0\in\thg$.

On the other hand, it follows from the above equation that for any
quotient $\tau_0$ of $\tau$, $\falbg{\tau_0}$
translation invariant, closed $*$-subalgebra of $\faltg$.
\endpf

To expand the next result, let us introduce another class of 
subspaces of $\fsalg$.  If $\tau\in\tg$ let
\[
\fsaltg=\fsal{G_\tau}\comp\eta_\tau.
\]
Just as with $\faltg$, we have that $\fsaltg\cong\fsal{G_\tau}$ completely
isometrically.

\begin{containment}\label{prop:containment}
If $\tau_1,\tau_2\in\tg$ the the following are equivalent:

{\rm (i)} $\tau_1\subseteq\tau_2$,

{\rm (ii)} $\falbg{\tau_1}\falbg{\tau_2}\subseteq\falbg{\tau_2}$,

{\rm (iii)} $\falbg{\tau_1\vee\tau_2}=\falbg{\tau_2}$,

{\rm (iv)} $\falbg{\tau_1}\subseteq\fsalbg{\tau_2}$, and

{\rm (v)} $\fsalbg{\tau_1}\subseteq\fsalbg{\tau_2}$.
\end{containment}

\proof (i)$\implies$(v) $\fsal{G_{\tau_2}}\supseteq
\fsal{G_{\tau_1}}\comp\eta^{\tau_2}_{\tau_1}$.  Hence
\[
\fsalbg{\tau_2}=\fsal{G_{\tau_2}}\comp\eta_{\tau_2}\supseteq
\fsal{G_{\tau_1}}\comp(\eta^{\tau_2}_{\tau_1}\comp\eta_{\tau_2})
=\fsal{G_{\tau_1}}\comp\eta_{\tau_1}=\fsalbg{\tau_1}.
\]

(v)$\implies$(iv) Trivial.

(iv)$\implies$(ii) $\falbg{\tau_2}$ is an ideal in $\fsalbg{\tau_2}$.

(ii)$\implies$(i) By Proposition \ref{prop:spinealgebra}
we obtain that $\falbg{\tau_1\vee\tau_2}\subseteq\falbg{\tau_2}$.
Hence we have that $\falbg{\tau_1\vee\tau_2}$ is a  
$*$-subalgebra of $\fal{G_{\tau_2}}$, which is translation invariant
modulo $N_\tau$.  
Thus, by Lemma \ref{lem:subalgebraq},
$\tau_1\vee\tau_2$ is a quotient of $\tau_2$, so
$\tau_1\vee\tau_2\subseteq\tau_2$.  This implies that $\tau_1\vee\tau_2
=\tau_2$, so $\tau_1\subseteq\tau_2$.

(i)$\implies$(iii) $\tau_1\vee\tau_2=\tau_2$.

(iii)$\implies$(ii) This follows Proposition \ref{prop:spinealgebra}. \endpf

\subsection{Decomposition of $\spineg$}

We cannot expect in the definition (\ref{eq:spinedef}) that a direct
sum decomposition obtains, as the next lemma shows.
We will see in Theorem \ref{theo:directsum}, that we may obtain a direct
sum if we omit non-trivial quotient topologies.

\begin{intersection}\label{lem:intersection}
If $\tau_1,\tau_2\in\tg$, then either $\falbg{\tau_1}\cap\falbg{\tau_2}=\{0\}$,
or each $\tau_j$ ($j=1,2$) is a quotient of a topology
$\taunq$ and
$\falbg{\tau_j}\subseteq\falbg{\taunq}$.
\end{intersection}

\proof Suppose $\fA=\falbg{\tau_1}\cap\falbg{\tau_2}\not=\{0\}$.
Then $\fA$ is a translation invariant, 
$*$-subalgebra of each $\falbg{\tau_j}$, $j=1,2$.
Thus it follows Lemma \ref{lem:subalgebraq} that $\fA=\falbg{\tau}$,
where $\tau$ is a mutual quotient of $\tau_1$ and $\tau_2$.
It thus follows Theorem \ref{theo:nonquot} (ii)
that each of $\tau_1$ and $\tau_2$ is a quotient of $\taunq$.  It thus 
follows Lemma \ref{lem:subalgebraq}, again, that 
$\falbg{\tau_j}\subseteq\falbg{\taunq}$ for each $j$.  \endpf

Recall the definition of $\tnqg$ is given in (\ref{eq:tnqg}).
We also define
\[
\tnqhg=\tnqg\cap\thg. 
\]

\begin{directsum}\label{theo:directsum}
The following direct sum decompositions obtain:
\[
\spineg=\ell^1\text{-}\!\!\!\!\bigoplus_{\tau\in\tnqg}\faltg
\quad\aand\quad
\rspineg=\ell^1\text{-}\!\!\!\!\bigoplus_{\tau\in\tnqhg}\faltg.
\]
\end{directsum}

\proof First note that by Theorem \ref{theo:nonquot}, for every
$\tau\iin\tg$ there is a $\taunq\iin\tnqg$ such that
$\tau$ is a quotient of $\taunq$.  It is clear that $\taunq\in\thg$
if $\tau\in\thg$.  Hence, by Lemma \ref{lem:subalgebraq},
$\faltg\subseteq\falbg{\taunq}$.  Thus $\wbar{\sum_{\tau\in\tnqg}\faltg}
=\spineg$ and $\wbar{\sum_{\tau\in\tnqhg}\faltg}=\rspineg$.  
These sums are direct sums, since by Lemma \ref{lem:intersection},
if $\tau_1,\tau_2\in\tnqg$, then $\falbg{\tau_1}\cap\falbg{\tau_2}=\{0\}$.
The direct sums are $\ell^1$-direct sums by \cite[3.13]{arsac}.  \endpf

We recall from \S \ref{subsec:latsemigrp} that $\tnqg$ and
$\tnqhg$ admit particular semigroup structures.

\begin{grading}\label{cor:grading}
{\bf (i)} $\spineg$ is a unital Banach algebra which
is {\rm graded} over the semilattice $\tnqg$, i.e.\ for
$\tau_1,\tau_2\iin\tnqg$, $\falbg{\tau_1}\falbg{\tau_2}\subseteq
\falbg{\tau_1\vee\tau_2}$.

{\bf (ii)} $\rspineg$ is graded over the semilattice  $\tnqhg$
and is an ideal in $\spineg$.  When $G$ is maximally almost periodic,
$\spineg=\rspineg$.
\end{grading}

\proof This follows from the theorem above, Lemma \ref{lem:subalgebraq}
Proposition \ref{prop:spinealgebra}, and \S \ref{subsec:latsemigrp}.  \endpf

It is shown in \cite{cowling} that for a semi-simple Lie group with finite
centre $G$, that $\fsalg$ can itself be written as a graded Banach algebra
over a finite lattice, whose summands have computable spectra.

\section{Characters and Homomorphisms}\label{sec:charhomo}

The general form of this section was inspired by the work of Inoue
\cite{inoue1,inoue2}.  All of the major results we obtain in the 
first 3 subsections below were obtained for abelian groups in \cite{inoue1}.  
Since any abelian group $G$
is maximally almost periodic, $\spineg=\rspineg$ in that case.
Despite that our results are more general than those of Inoue, our proofs
are often simplified by our context.

\subsection{The Spectrum of $\spineg$} \label{ssec:spectrum}
We recall that $\tnqg$, as defined in 
(\ref{eq:tnqg}), is a directed set via inclusion.  We say that
a subset $\fS$ of $\tnqg$ is {\it hereditary} if for any $\tau_0$
in $\fS$, the set
\begin{equation}\label{eq:principalhd}
\fS_{\tau_0}=\{\tau\in\tnqg:\tau\subseteq\tau_0\}
\end{equation}
is contained in $\fS$.  Notice that each set $\fS_{\tau_0}$
is hereditary and directed.  We define
\[
\hdg=\{\fS\subseteq\tnqg:\fS\text{ is hereditary, directed and non-empty}\}.
\]
Now if $\fS\in\hdg$, we let $G_\fS={\underleftarrow{\lim}}_{\tau\in\fS}G_\tau$
as in (\ref{eq:projlim}).  We then let
\[
G^*=\bigsqcup_{\fS\in\hdg}G_\fS
\]
where we use $\sqcup$ to denote coproduct. 

We have from Theorem \ref{theo:directsum} that if $u\in\spineg$, then
there exists for each $\tau\iin\tnqg$ a unique $u_\tau\iin\faltg$ such that
\begin{equation}\label{eq:indsumdecomp}
u=\sum_{\tau\in\tnqg}u_\tau\quad\aand\quad\norm{u}=
\sum_{\tau\in\tnqg}\norm{u_\tau}.
\end{equation}
For each $\tau$, (\ref{eq:faltg}) shows 
that there exists a unique $\hat{u}_\tau\iin\fal{G_\tau}$ such that
$\hat{u}_\tau\comp\eta_\tau=u_\tau$.  In fact, $u_\tau\mapsto
\hat{u}_\tau$ is the Gelfand transform of $\faltg$.

A {\it character} on an abelian complex algebra
 is a non-zero multiplicative linear functional
into the complex number field $\Cee$.

\begin{spectrum}\label{theo:spectrum}
{\bf (i)} If $s\in G^*$, say $s=(s_\tau)_{\tau\in\fS}\in G_\fS$ for
some $\fS$ in $\hdg$, then the functional defined by
\[
\chi_s\left(\sum_{\tau\in\tnqg}u_\tau\right)
=\sum_{\tau\in\fS}\hat{u}_\tau(s_\tau)
\]
is a character on $\spineg$.

{\bf (ii)} If $\chi:\spineg\to\Cee$ is a character, then there exists
an $s\iin G^*$ such that $\chi=\chi_s$.

Thus the map $s\mapsto\chi_s$ is a bijection from $G^*$ into the Gelfand
spectrum of $\spineg$.
\end{spectrum}

\proof  Our proof will rely on the following formula:
if $\tau_1,\tau_2\in\tnqg$, $u_j\in\falbg{\tau_j}$ for $j=1,2$,
and $t\in G_{\tau_1\vee\tau_2}$ then
\begin{equation}\label{eq:basicformula}
\what{u_1u_2}(t)=
\hat{u}_1\bigl(\eta^{\tau_1\vee\tau_2}_{\tau_1}(t)\bigr)
\hat{u}_2\bigl(\eta^{\tau_1\vee\tau_2}_{\tau_2}(t)\bigr).
\end{equation}
To see this let $(t_i)$ be a net from $G$ such that 
$\lim_i \eta_{\tau_1\vee\tau_2}(t_i)=t$.  Then for
each $j=1,2$ we have that $\lim_i \eta_{\tau_j}(t_i)
=\lim_i \eta^{\tau_1\vee\tau_2}_{\tau_j}\comp
\eta_{\tau_1\vee\tau_2}(t_i)=
\eta^{\tau_1\vee\tau_2}_{\tau_j}(t)$.  Hence
\begin{align*}
\what{u_1u_2}(t) &=
\lim_i \what{u_1u_2}\bigl(\eta_{\tau_1\vee\tau_2}(t_i)\bigr)
=\lim_i u_1u_2(t_i) 
=\lim_i u_1(t_i)u_2(t_i) \\
&=\lim_i \hat{u}_1\bigl(\eta_{\tau_1}(t_i)\bigr)
\hat{u}_2\bigl(\eta_{\tau_2}(t_i)\bigr) 
=\hat{u}_1\bigl(\eta^{\tau_1\vee\tau_2}_{\tau_1}(t)\bigr)
\hat{u}_2\bigl(\eta^{\tau_1\vee\tau_2}_{\tau_2}(t)\bigr).
\end{align*}
We note that if $\tau_1\subseteq\tau_2$, then (\ref{eq:basicformula})
reduces to
\begin{equation}\label{eq:basicformula1}
\what{u_1u_2}(t)=\hat{u}_1\bigl(\eta^{\tau_2}_{\tau_1}(t)\bigr)\hat{u}_2(t).
\end{equation}

(i) Since $\fS$ is non-empty, $\chi_s|_{\faltg}\not=0$ for some
$\tau\iin\fS$, so $\chi_s\not=0$.  Using (\ref{eq:indsumdecomp})
it is straightforward to verify that $\chi_s$ is linear and continuous.
It remains to verify that $\chi_s$
is multiplicative.  If $u,v\in\spineg$, then by (\ref{eq:indsumdecomp})
we obtain
\[
uv=\left(\sum_{\tau_1\in\tnqg}u_{\tau_1}\right)
\left(\sum_{\tau_2\in\tnqg}v_{\tau_2}\right)
=\sum_{\tau_1\in\tnqg}\sum_{\tau_2\in\tnqg}u_{\tau_1}v_{\tau_2}
\]
where the double sum converges absolutely.  Since by Corollary 
\ref{cor:grading}, $u_{\tau_1}v_{\tau_2}\in\falbg{\tau_1\vee\tau_2}$
we obtain that
\[
\chi_s(uv)=\sum_{\tau_1\in\tnqg}\sum_{\tau_2\in\tnqg}
\what{u_{\tau_1}v_{\tau_2}}(s_{\tau_1\vee\tau_2}).
\]
By definition of $G_\fS$ we see that 
$\eta^{\tau_1\vee\tau_2}_{\tau_j}(s_{\tau_1\vee\tau_2})
=s_{\tau_j}$ for $j=1,2$.  It then follows from (\ref{eq:basicformula})
that
\[
\chi_s(uv)=\sum_{\tau_1\in\tnqg}\sum_{\tau_2\in\tnqg}
\hat{u}_{\tau_1}(s_{\tau_1})\hat{v}_{\tau_2}(s_{\tau_2})
\]
which is easily seen to be $\chi_s(u)\chi_s(v)$.

Any pair of distinct points $s$ and $s'$
in $G^*$ form distinct characters $\chi_s$ and $\chi_{s'}$.
Indeed, if $s_\tau\not=s'_\tau$, then
$\chi_s|_{\faltg}\not=\chi_{s'}|_{\faltg}$.

(ii)  If $\chi:\spineg\to\Cee$ is a character then let
\[
\fS_\chi=\left\{\tau\in\tnqg:\chi|_{\faltg}\not=0\right\}.
\]
It follows from Corollary \ref{cor:grading} that $\fS_\chi$ is hereditary
and directed.
Since for each $\tau\iin\fS_\chi$, $\chi|_{\faltg}$ is a character,
we may define $s_\tau$ by the relation
\[
\hat{u}(s_\tau)=\chi(u)
\]
for each $u\iin\faltg$.  If $\tau_1\subseteq\tau_2$ in $\fS_\chi$ then
we have for $u_j\in\falbg{\tau_j}$, $j=1,2$, that
\[
\hat{u}_1(s_{\tau_1})\hat{u}_2(s_{\tau_2})
=\chi(u_1)\chi(u_2)=\chi(u_1u_2)=\what{u_1u_2}(s_{\tau_2}).
\]
Applying (\ref{eq:basicformula1}), we then obtain
\[
\hat{u}_1(s_{\tau_1})\hat{u}_2(s_{\tau_2})=
\hat{u}_1\bigl(\eta^{\tau_2}_{\tau_1}(s_{\tau_2})\bigr)\hat{u}_2(s_{\tau_2}).
\]
Since this holds for all choices $u_j\in\falbg{\tau_j}$, $j=1,2$, it follows
that $s_{\tau_1}=\eta^{\tau_2}_{\tau_1}(s_{\tau_2})$.  Hence
$s=(s_\tau)_{\tau\in\fS_\chi}\in G_{\fS_\chi}$, and $\chi=\chi_s$.  \endpf

\begin{spectopology}\label{theo:spectopology}
Let a topology on $G^*$ be given as follows:  given any $s_0\iin G^*$,
say $s_0\in G_{\fS_0}$ for some $\fS_0\iin\hdg$, a neighbourhood
basis at $s_0$ is formed by
\[
\left\{U(V_{\tau_0};W_{\tau_1},\dots,W_{\tau_n}):
\begin{matrix} \tau_0\in\fS_0, V_{\tau_0}\text{ is an open nbhd.\ of }
s_{0,\tau_0}\iin G_{\tau_0}; \\
\tau_j\notin\fS_0, G_{\tau_j}\setdif W_{\tau_j}\text{ is compact in }
G_{\tau_j},j=1,\dots,n \end{matrix}\right\}
\]
where
\[
U(V_{\tau_0};W_{\tau_1},\dots,W_{\tau_n})=\left\{s\in G^*:
\begin{matrix}
s\in G_\fS\text{ for some }\fS\supseteq\fS_0\iin\hdg, \\
\text{ for which }s_{\tau_0}\in V_{\tau_0},\text{ and } \\
s_{\tau_j}\in  W_{\tau_j}\text{ if }\tau_j\in\fS,
\text{ for }j=1,\dots,n
\end{matrix}\right\}.
\]
Then this is the topology $G^*$ inherits as the spectrum of $\spineg$.
\end{spectopology}

\proof We note that for a net $(s_i)\iin G^*$ the following are equivalent:

(i) $s_i\to s_0$ in $G^*$ with topology as described above,

(ii) $\displaystyle s_{i,\tau}\to
\begin{cases} s_{0,\tau} &\text{if }\tau\in\fS_0, \\
\infty &\text{if }\tau\not\in\fS_0 \end{cases}$ and

(iii) $\chi_{s_i}\to\chi_{s_0}$ weak* in $\spineg^*$.

\noindent
The equivalence of (i) and (ii) is clear.  
Now, let $u\in\spineg$ and write $u=\sum_{\tau\in\tnqg}u_\tau$ as in
(\ref{eq:indsumdecomp}). Then we  see that
\[
\chi_{s_i}(u)=\left(\sum_{\tau\in\fS_0}\hat{u}_\tau(s_{i,\tau})
+\sum_{\tau\in\fS_i\setdif\fS_0}\hat{u}_\tau(s_{i,\tau})\right)\overset{i}{\to}
\sum_{\tau\in\fS_0}\hat{u}_\tau(s_{0,\tau})=\chi_{s_0}(u)
\]
where $s_i\in G_{\fS_i}$ for each $i$.
This shows that (iii) follows (ii).   The same calculation shows that
(iii) implies (ii), if we select $u=u_\tau$ for a single $\tau$.  \endpf

\begin{spectopology1}\label{cor:spectopology1} Let $\fS\in\hdg$.  Then

{\bf (i)} the injection $G_{\fS}\hookrightarrow G^*$ is continuous when 
$G_{\fS}$ has the projective limit topology $\tau_\fS
=\bigvee_{\tau\in\fS}\tau$, and

{\bf (ii)} $\displaystyle \wbar{G_\fS}
=\bigsqcup_{\fS_0\in\hdg,\fS_0\subseteq\fS}G_{\fS_0}$.

\noindent In particular, the injection $G\to G^*:
s\mapsto\bigl(\eta_\tau(s)\bigr)_{\tau\in\tnqg}$, which is the realisation of 
$G$ as evaluation functionals, 
has dense range.
\end{spectopology1}

\proof (i) Clear, by definition.  (ii)  Any of the open neighbourhoods
of a point $s_0\iin\fS_0\subseteq\fS$ described in the above theorem
clearly intersect $G_\fS$.  Conversely, if $s_0\in G_{\fS_0}$ where
$\fS_0\not\subset\fS$, and $\tau_0\in\fS_0\setdif\fS$, then
for any neighbourhood $V_{\tau_0}$ of $s_{\tau_0}$ in $G_{\tau_0}$
\[
U(V_{\tau_0})=\left\{s\in G^*:s\in G_{\fS'}\text{ for some }\fS'\supset\fS_0
\aand s_{\tau_0}\in V_{\tau_0}\right\}
\]
is a neighbourhood of $s_0$ which misses $G_\fS$.
\endpf

Let us observe that $G^*$ is a semi-topological semigroup.  
Let
\[
\lam_*=\bigoplus_{\tau\in\tnqg}\lam_\tau:
G\to\fB\left({\ell^2\text{-}\bigoplus}_{\tau\in\tnqg}\bltwo{G_\tau}\right).
\]
Then $\ccfs{{\lam_*}}=\spineg$, so by \cite[2.2]{arsac},
the von Neumann algebra $\vn{{\lam_*}}=\lam_*(G)''$
is the dual of $\spineg$.  It follows from Corollary \ref{cor:spectopology1}
that we obtain a bicontinuous semigroup
isomorphism $\wbar{\lam_*(G)}^{w^*}\cong G^*$ which continuously
extends the map $\lam_*(s)\mapsto\bigl(\eta_\tau(s)\bigr)_{\tau\in\tnqhg}$.  
The unit ball of $\fB(\fH)$, where $\fH$ is a Hilbert space,
is a semi-topological semigroup in the $w^*$ topology.  Hence the 
$w^*$-closure of any sub(semi)group is a semi-topological
semigroup.  Let us thus call $G^*$ the
{\it spine compactification} of $G$.
We note that the multiplication on $G^*$ is given by
\begin{equation}\label{eq:gstarmult}
(s_\tau)_{\tau\in \fS_1}(t_\tau)_{\tau\in \fS_2}
=(s_\tau t_\tau)_{\tau\in \fS_1\cap\fS_2}
\end{equation}
where $\fS_1,\fS_2\in\hdg$.  
Hence we have that $G_{\fS_1}G_{\fS_2}
= G_{\fS_1\cap\fS_2}$.

\subsection{$\hdg$ is Dual to $\tnqg$}\label{subsec:dual}
The set $\hdg$ is a semilattice under the 
operation $(\fS_1,\fS_2)\mapsto\fS_1\cap\fS_2$.  We note that
these intersections are always non-empty, for it follows from
Theorem \ref{theo:nonquot}
that $\tauap\in\fS$ for every $\fS\in\hdg$.

The semilattice $\hdg$ is dual
to the semilattice $\tnqg$.  Indeed, if $\omega:\tnqg\to\Cee$
is a non-zero semicharacter -- i.e.\ $\omega(\tau_1\vee\tau_2)
=\omega(\tau_1)\omega(\tau_2)$ -- then, since $\tnqg$ is an idempotent 
semigroup, we must have that $\omega=1_{\fS_\omega}$, the indicator function of
a nonempty subset $\fS_\omega$ of $\tnqg$.  It is straightforward to check that
$\fS_\omega$ is hereditary and that if $\tau_1,\tau_2\in\fS_\ome$, then
$\tau_1\vee\tau_2\in\fS_\ome$ too.  Hence $\fS_\ome\in\hdg$.  Moreover,
the dual pairing $\dpair{\cdot}{\cdot}:\tnqg\cross
\hdg\to\{0,1\}$ given by
\[
\dpair{\tau}{\fS}=1_\fS(\tau)=\begin{cases} 1 &\text{if }\tau\in\fS \\
0 &\text{if }\tau\not\in\fS \end{cases}
\]
satisfies
\[
\dpair{\tau_1\vee\tau_2}{\fS}=
\dpair{\tau_1}{\fS}\dpair{\tau_2}{\fS}\quad\aand\quad
\dpair{\tau}{\fS_1\cap\fS_2}=\dpair{\tau}{\fS_1}\dpair{\tau}{\fS_2}
\]
for $\tau_1,\tau_2,\tau\iin\tnqg$ and $\fS_1,\fS_2,\fS\iin\hdg$.

On the other hand, if $\Omega:\hdg\to\Cee$ is a non-zero semicharacter
then $\Omega=1_{\mathfrak{S}}$ for some $\mathfrak{S}\subseteq\hdg$.
The set $\fS_\Omega=\bigcap_{\fS\in\mathfrak{S}}\fS$ is non-empty and thus 
itself in $\hdg$.  Thus we have
\[
\Omega(\fS)=\begin{cases} 1 &\text{if }\fS\supseteq\fS_\Omega \\
0 &\text{if }\fS\not\supseteq\fS_\Omega \end{cases}.
\]
We may express this as ``$\Omega(\fS)=\dpair{\tau_{\fS_\Omega}}{\fS}$'',
where $\tau_{\fS_\Omega}=\bigvee_{\tau\in\fS_\Omega}\tau$.
Thus, the lattice completion of $\tnqg$, 
$\ctnqg=\{\tau_\fS:\fS\in\hdg\}$, is the dual of $\hdg$.

Since any semilattice is a commutative {\it inverse semigroup},
it is possible to equip $\hdg$, qua the dual of $\tnqg$, with a compact
semigroup topology -- the topology of pointwise convergence as characters --
such that $\tnqg$ represents the set of continuous
characters on $\hdg$.  In other words, Pontryagin duality holds.
See \cite{austin} for details.  In analogy with
Proposition \ref{prop:idempotents}, below, we can see that $\hdg$
is semilattice isomorphic to the idempotent lattice $E$ in $G^*$,
which is a compact subsemigroup.
It is easy to verify that the compact semilattice $\hdg$ is homeomorphic
to $E$. It is curious that in many situations,
as in Proposition \ref{prop:hdprincipal} below, that all semicharacters
on $\hdg$ are automatically continuous.

We note that the semilattice structure on $\hdg$ gives rise to
an infemum construction on $\ctnqg$:
\[
\tau_{\fS_1}\wedge\tau_{\fS_2}=\tau_{\fS_1\cap\fS_2}.
\]
Thus if $\tau_1,\tau_2\in\tnqg$, then we then have
\begin{equation}\label{eq:inftopology}
\tau_1\wedge\tau_2=\tau_{\fS_{\tau_1}\cap\fS_{\tau_2}}
\end{equation}
where $\fS_{\tau_j}$ ($j=1,2$) are the principal hereditary directed sets
from (\ref{eq:principalhd}).

\begin{hdprincipal}\label{prop:hdprincipal}
The following are equivalent:

{\bf (i)} every $\fS\iin\hdg$ is principal, and

{\bf (ii)} $\tnqg$ is a complete semilattice.

\noindent  In this case
$\tnqg$ is a lattice, i.e.\ if $\tau_1,\tau_2\in\tnqg$,
then $\tau_1\wedge\tau_2\in\tnqg$ too.
Moreover, the semilattice $(\tnqg,\wedge)$ is lattice 
isomorphic to $\hdg$.
\end{hdprincipal}

It is shown in \S \ref{ssec:freegrp} that condition (ii) does not always hold.

\medskip
\proof (i) $\implies$ (ii)  If $\fF$ is any family from $\tnqg$,
let $\langle\fF\rangle=\bigcup\{
\fS_{\tau_1\vee\dots\vee\tau_n}:\tau_1,\dots,\tau_n\in\fF, n=1,2,\dots\}$
where each $\fS_{\tau_1\vee\dots\vee\tau_n}$ is the
principal set from (\ref{eq:principalhd}).  Then 
$\langle\fF\rangle\in\hdg$ and hence there is $\tau_0$ such that
$\langle\fF\rangle=\fS_{\tau_0}$.  Thus
\[
\bigvee_{\tau\in\fF}\tau=\bigvee_{\tau\in\langle\fF\rangle}\tau
=\tau_{\fS_{\tau_0}}=\tau_0\in\tnqg.
\]
The converse, (ii) $\implies$ (i), is trivial.

We obtain the infemum operation on $\tnqg$ by virtue of
(\ref{eq:inftopology}).  The lattice isomorphism
$(\tnqg,\wedge)\cong(\hdg,\cap)$ follows as well.  \endpf

Thus if $\tnqg$ is a complete lattice we have
\begin{equation}\label{eq:tnqgcompletespec}
G^*=\bigsqcup_{\tau\in\tnqg}G_\tau.
\end{equation}
In this case, if $s\in G^*$, so $s\in G_{\tau_0}$ for some $\tau_0\iin\tnqg$,
the character $\chi_s$ it implements on $\spineg$, as in
Theorem \ref{theo:spectrum} (i), is
\[
\chi_s(u)=\sum_{\tau\in\fS_{\tau_0}}
\hat{u}_\tau\bigl(\eta^{\tau_0}_\tau(s)\bigr).
\]
Moreover, the semigroup operation (\ref{eq:gstarmult}) admits the
following formula:  if $s_j\in G_{\tau_j}$, $\tau_j\in\tnqg$, $j=1,2$, then
\begin{equation}\label{eq:gstarmult1}
s_1s_2=\eta^{\tau_1}_{\tau_1\wedge\tau_2}(s_1)
\eta^{\tau_2}_{\tau_1\wedge\tau_2}(s_2).
\end{equation}
Thus we have $G_{\tau_1}G_{\tau_2}\subseteq G_{\tau_1\wedge\tau_2}$.
In effect, $G^*$ is graded over the semilattice $(\tnqg,\wedge)$.

\subsection{The Spectrum of $\rspineg$}\label{ssec:rspinespec}
Let $\fS$ be a directed subset of $\tnqhg$.  We say that
$\fS$ is {\it hereditary in $\tnqhg$} if $\fS=\fS'\cap\tnqhg$
for some $\fS'\iin\hdg$.  We then let
\[
\rhdg=\{\fS\cap\tnqhg:\fS\in\hdg\} 
\]
We note that if $G$ in not maximally almost periodic, then
$\varnothing=\{\tauap\}\cap\tnqhg\in\rhdg$.
Just as in \S\ref{subsec:dual}, above, we
find that $\rhdg$ is the dual semilattice to $\tnqhg$.
We now let
\begin{equation}\label{eq:rspinespec}
G^*_0=\bigsqcup_{\fS\in\rhdg}G_\fS
\end{equation}
where $G_\fS$ is as in (\ref{eq:projlim}) if $\fS\not=\varnothing$
and $G_\varnothing=\{\infty\}$.  

The following is proved just as in Theorem \ref{theo:spectrum}.

\begin{rspectrum}\label{theo:rspectrum}
There is a bijective correspondence between

{\bf (i)} $G^*_0$, and

{\bf (ii)} the set of multiplicative linear functionals on $\rspineg$,

\noindent as given in Theorem \ref{theo:spectrum}.  If 
$\varnothing\in\rhdg$, then $\infty$ corresponds to the zero functional.
\end{rspectrum} 

We gain a description of the compact topology on $G^*_0$ inherits
from being the spectrum (or the spectrum's one point compactification) 
which is similar to Theorem \ref{theo:spectopology}.  
We also obtain a multiplication on $G^*_0$ as 
we do for $G^*$:  if 
$s_j=(s_{j,\tau})_{\tau\in\fS_j}\in G_{\fS_j}$, $j=1,2$, then
\[
s_1s_2=(s_{1,\tau}s_{2,\tau})_{\tau\in\fS_1\cap\fS_2}
\]
where we write $s_1s_2=\infty$ if $\fS_1\cap\fS_2=\varnothing$.
We thus call $G^*_0$ the {\it reduced spine compactification} of $G$.

\begin{idempotents}\label{prop:idempotents}
{\bf (i)} Every idempotent in $G^*_0$ is of the form 
$e_\fS=\bigl(\eta_\tau(e)\bigr)_{\tau\in\fS}$ for some $\fS\in\rhdg$,
where $e_\varnothing=\infty$.

{\bf (ii)} The set of idempotents in $G^*_0$, $E_0$ forms a semilattice
which is isomorphic to the semilattice $\rhdg$.  Moreover,
$E_0$ is a closed subsemigroup and thus a topological semilattice.

{\bf (iii)} $\displaystyle e_\fS G^*_0=G^*_0 e_\fS=\bigsqcup_{\fS_0\in\rhdg,
\fS_0\subseteq\fS}G_{\fS_0}=\wbar{G_\fS}$, the closure of $G_\fS$
in $G^*_0$.
\end{idempotents}

\proof (i) If $f$ is an idempotent in $G^*$, then $f=(s_\tau)_{\tau\in\fS}$
for some $\fS\in\hdg$.  But then $s_\tau^2=s_\tau$ for each $\tau$, so
$s_\tau=\eta_\tau(e)$. 

(ii) The fact that $E_0$ is a semilattice, isomorphic to $\rhdg$, is clear.  
To see that it is closed, let $(e_i)$ be a net in $E_0$ 
converging to $s\iin G^*_0$, so $s=(s_\tau)_{\tau\in\fS}$ for some $\fS$. 
Then we have from the $G^*_0$-analogue of Theorem \ref{theo:spectopology}
that $e_i\in G_\fS$ for large $i$ and $e_{i,\tau}\to s_\tau$ for each 
$\tau\in\fS$.  It then follows that $s=e_\fS$, as in (i).
That $E_0$ is a semi-topological semigroup then follows from 
\cite[1.4.11]{berglundjm}.

(iii)  The first two equalities
are obvious.  The third is that same as 
Corollary \ref{cor:spectopology1} (ii).
\endpf

These idempotents correspond in an obvious way
to certain central projections in $\vn{{\lam_*^0}}=
\ell^\infty\text{-}\bigoplus_{\tau\in\tnqhg}\vn{{\lam_{\tau}}}$.
It thus follows that $G^*_0$ is a {\it Clifford semigroup} presented
in a standard from (see \cite[Theorem 3]{clifford} for an intrinsic
description mirroring (\ref{eq:rspinespec}), or 
\cite[Proposition 2.1.4]{pattersonG} for a representation theorem
of general inverse semigroups).
The idempotents of $G^*$ can be obtained similarly, and it too
is a Clifford semigroup.

Recall from Corollary \ref{cor:grading}
that $\rspineg$ is unital if $G$ is maximally almost periodic.
This characterises maximal almost periodicity.

\begin{whenunital}\label{theo:whenunital}
If $\rspineg$ is unital then $G$ is maximally almost periodic.
\end{whenunital}

\proof If $\rspineg$ is unital, then 
by Theorem \ref{theo:rspectrum} $G^*_0$ is compact
and does not contain the zero functional.
By Proposition \ref{prop:idempotents} (ii) the family $E_0$ of 
idempotents is a compact semilattice, and hence 
there is a minimal idempotent, $e_m$.  Indeed, $e_m$ is a cluster
point of the net $(e)_{e\in E_0}$, partially ordered by $e\leq e'$
if and only if $ee'=e'$; this is the usual partial order, reversed. 
It also follows from Proposition \ref{prop:idempotents} (i)
that $e_m=e_{\fS_m}$ where $\fS_m\in\rhdg$.  But, by
minimality, $\fS_m$ is a singleton $\{\tau_m\}$.  Then by
Proposition \ref{prop:idempotents} (iii), $e_mG^*_0=G_{\tau_m}$ is a compact
group, whence $G$ embeds in a compact group.  Since $\tau_m\in\tnqg$,
$\tau_m=\tauap$.  \endpf

\subsection{Homomorphisms into $\fsalh$}
\label{ssec:homomorphisms}
The major result of this section is a generalisation of
the present authors' theorem characterising completely bounded
homomorphisms from $\falg$ (when $G$ is amenable) to $\fsalh$.
This result generalised the major result of \cite{cohen}.
The homomorphism theorem we obtain for spines generalises
that of \cite{inoue1}.

A full account of piecewise affine maps between groups in given in
\cite{ilies}, or in \cite{rudin} in the abelian case.  We summarise
some basic facts below.

Let $H$ and $G_0$ be groups.  Let $\cringh$ denote the {\it coset ring}
of $H$, that is, the smallest ring of subsets containing all cosets of
subgroups of $H$.  A map $\alp:Y\subseteq H\to G_0$ is called
{\it piecewise affine} if

\medskip
{\bf (i)} there are pairwise disjoint $Y_j\iin\cringh$, $j=1,\dots,n$
such that $Y=\dot{\bigcup}_{j=1}^n Y_j$, and

{\bf (ii)} each $Y_j$ is contained in a coset $L_j$ on which there
is an affine map $\alp_j:L_j\to G_0$ such that $\alp_j|_{Y_j}=\alp|_{Y_j}$.

\medskip
\noindent We recall that an {\it affine} map $\alp:C\subseteq H\to G_0$,
where $C$ is a coset, is a map which satisfies
\[
\alp(rs^{-1}t)=\alp(r)\alp(s)^{-1}\alp(t)
\]
for any $r,s,t\iin C$.  We recall too that cosets are characterised by the
condition: $r,s,t\iin C$ implies $rs^{-1}t\in C$.

The following is proved exactly as Lemma 1.3 (ii) in \cite{ilies}.
If $H$ is a topological group we will denote by $\ocringh$ the
smallest ring of subsets containing all open cosets.

\begin{pamaps}\label{lem:pamaps}
Let $H$ be a locally compact group and $G_0$ be a complete Hausdorff 
topological group.  The any continuous piecewise affine map 
$\alp:Y\subseteq H\to G_0$ admits a continuous extension
$\bar{\alp}:\wbar{Y}\subseteq H\to G_0$.  Moreover, $\wbar{Y}$ is open
and admits a decomposition as in (ii) above,  
$\wbar{Y}=\dot{\bigcup}_{j=1}^n Y_j$, where each $Y_j\in\ocringh$.
\end{pamaps}

Thus we may hereafter assume that all continuous piecewise affine maps
between such pairs of groups have open and closed domains.  If
$\alp=\bar{\alp}$ is as above, let us write
\begin{equation}\label{eq:pabound}
c(\alp,Y)=n\left(\sum_{j=1}^n\norm{1_{Y_j}}\right)
\end{equation}
where $\norm{\cdot}$ is the norm on $\fsalh$.

We will now assume for the rest of the section that $H$ is a locally compact 
group.  A map $\alp:Y\subseteq H\to G^*$ is called {\it piecewise affine}
if

{\bf (i)} there is an $\fS\iin\hdg$ such that $\alp(Y)\subseteq G_\fS$, and

{\bf (ii)} $\alp:Y\subseteq H\to G_\fS$ is continuous and piecewise affine.

\noindent
Recall that $G_\fS={\underleftarrow{\lim}}_{\tau\in\fS}G_\tau$ is complete,
so the lemma above applies.

\begin{homomorphism}\label{theo:homomorphism}
{\bf (i)} If $\alp:Y\subseteq H\to G^*$ is piecewise affine then
$\Psi_\alp:\spineg\to\fsalh$, given for each $u\iin\spineg$ and $h\iin H$ by
\[
\Psi_\alp u(h)=\begin{cases} \chi_{\alp(h)}u &\text{if }h\in Y, \\
0 &\text{otherwise} \end{cases}
\]
is a completely bounded homomorphism.

{\bf (ii)} If $G$ is amenable, then any completely bounded homomorphism
$\Psi:\spineg\to\fsalh$ is of the form $\Psi=\Psi_\alp$ for
some piecewise affine map $\alp:Y\subseteq H\to G^*$.
\end{homomorphism}

\proof (i) Say $\alp(H)\subseteq G_\fS$.  For each $\tau\iin\fS$ we
let $\eta^\fS_\tau:G_\fS\to G_\tau$ be the continuous injective homomorphism
which is guaranteed by the projective limit construction; hence
if $\tau_1\subseteq\tau_2\iin\fS$, then $\eta^{\tau_2}_{\tau_1}\comp
\eta^\fS_{\tau_2}=\eta^\fS_{\tau_1}$. Then
$\alp_\tau=\eta^\fS_\tau\comp\alp:Y\to G_\tau$ is continuous and piecewise
affine.  Moreover, we obtain that if $\tau_1\subseteq\tau_2\iin\fS$, then
\begin{equation}\label{eq:alptau}
\eta^{\tau_2}_{\tau_1}\comp\alp_{\tau_2}=\alp_{\tau_1}.
\end{equation}
Then for $u_\tau\iin\faltg$ and $h\iin H$ we have
\[
\Psi_\alp u_\tau(h)=\begin{cases} \hat{u}_\tau\bigl(\alp_\tau(h)\bigr)
&\text{if }h\in Y\aand \tau\in\fS, \\
0 &\text{otherwise} \end{cases}.
\]
Hence by \cite{ilies}, Proposition 3.1, $\Psi_\alp|_{\faltg}$ is a completely
bounded homomorphism with $\cbnorm{\Psi_\alp|_{\faltg}}\leq c(\alp,Y)$,
where $c(\alp,Y)$ is defined in (\ref{eq:pabound}).  It then follows
Theorem \ref{theo:directsum} that $\Psi_\alp$ extends to a linear
completely bounded map on  all of $\spineg$.  Finally, 
it follows from (\ref{eq:alptau}) and (\ref{eq:basicformula})
that  $\Psi_\alp$ is a homomorphism on $\spineg$.

%


%
(ii) Let $\fS_\Psi=\{\tau\in\tnqg:\Psi|_{\faltg}\not=0\}$.  It follows
from Corollary \ref{cor:grading} that $\fS_\Psi$ is hereditary
and directed.  Hence if $\Psi\not=0$, then $\fS_\Psi\in\hdg$.
We note that for any $\tau\iin\tg$, $G_\tau$ is amenable by
\cite[Proposition 1.2.1]{rundeB}.
Thus, by the main result of \cite{ilies}, if $\tau\in\fS_\Psi$ then
$\Psi|_{\faltg}$ takes the form
\[
\Psi u_\tau (h)=\begin{cases} \hat{u}_\tau\bigl(\alp_\tau(h)\bigr)
& \text{if }h\in Y_\tau, \\ 0 &\text{otherwise}\end{cases}
\]
for each $u_\tau\iin\faltg$, $h\iin H$ where 
$\alp_\tau:Y_\tau\subseteq H\to G_\tau$ is a continuous piecewise affine
map.  Since $\Psi$ is a homomorphism, it follows (\ref{eq:basicformula1})
that the system $\{\alp_\tau: Y_\tau\to G_\tau\}_{\tau\in\fS_\Psi}$ satisfies
\[
\text{if }\tau_1\subseteq\tau_2\iin\fS_\Psi\text{ then }
\eta^{\tau_2}_{\tau_1}\comp\alp_{\tau_2}=\alp_{\tau_1}.
\]
Moreover, the domains $Y_{\tau_j}$ ($j=1,2$) must then be equal.
Thus we can unambiguously denote each $Y_\tau$ by $Y$.
Thus there is a map $\alp:Y\to G_{\fS_\Psi}$ given by 
\[
\alp(h)=
\bigl(\alp_\tau(h)\bigr)_{\tau\in\fS_\Psi}.  
\]
This map is also piecewise
affine.  Indeed, suppose $\tau_1\subseteq\tau_2$ in $\fS_\Psi$,
$Y=Y_{\tau_1}=\dot{\bigcup}_{j=1}^nY_j$, and there exist for each $j=1,\dots,n$
open cosets $L_j\subseteq Y_j$ and affine maps $\alp_{\tau_1,j}:L_j\to
G_{\tau_1}$ extending $\alp_{\tau_1}|_{Y_j}$.  Then 
$\alp_{\tau_2,j}=\eta^{\tau_2}_{\tau_1}\comp\alp_{\tau_1,j}$ will be an
affine extension of $\alp_{\tau_2}|_{Y_j}$.  \endpf

\subsection{The Range of Homomorphisms}

The goal of this section is to establish that for any
locally compact group $H$, $\spineh$ is that part of $\fsal{H}$
which contains the ranges of
each of a natural class of homomorphisms from $\falg$, where
$G$ is another locally compact group.  For abelian groups
the was proved in \cite{inoue2}.

Much of the work from this section is done in the following technical lemma.

\begin{exttop}\label{lem:exttop}
Let $H_0$ be an open subgroup of a locally compact group $H$
and $\tau\in\tho$.  Then there is a $\bar{\tau}\iin\th$
such that $\bar{\tau}\cap H_0=\{U\cap H_0:U\in\bar{\tau}\}=\tau$,
and $(H_0)_\tau$ is bicontinuously isomorphic to 
$\wbar{\eta_{\bar{\tau}}(H_0)}$, which is an open subgroup
of $H_{\bar{\tau}}$.  Moreover, $\bar{\tau}\in\thh$ if $\tau\in\thho$.
\end{exttop}

\proof Let
\[
\bar{\tau}=\bigl\{U\subseteq H:s^{-1}(U\cap sH_0)\in\tau
\text{ for every }s\iin G\bigr\}.
\]
Then $\bar{\tau}\cap H_0=\tau$, and $\bar{\tau}$ is a locally
precompact group topology on $H$, since it has a neighbourhood basis 
of precompact neighbourhoods about $e$ in $H_0$.  It is clear that
$\wbar{\eta_{\bar{\tau}}(H_0)}\cong (H_0)_\tau$.
It is trivial to note that $\bar{\tau}$ is Hausdorff if $\tau$ is.

Now suppose that $s_0\in\wbar{\eta_{\bar{\tau}}(H_0)}$.  Let
$(t_j)$ be a net from $H$ such that $\lim_j\eta_{\bar{\tau}}(t_j)=s_0$.
Let $(s_i)$ be a net from $H_0$ such that 
$\lim_i \eta_\tau(s_i)=s_0$.  Then if we consider the product
net $(s_i^{-1}t_j)$ we have that $\lim_{(i,j)}\eta_{\bar{\tau}}(s_i^{-1}t_j)
=e_\tau$, whence $\bar{\tau}\text{-}\lim_{(i,j)}s_i^{-1}t_j=e$.
Thus we see that for large $(i,j)$ we have
$s_i^{-1}t_j\in H_0$, so for large $j$ we obtain $t_j\in s_iH_0=H_0$.
Hence the net $(t_j)$ is eventually contained in $H_0$, so
$\bigl(\eta_{\bar{\tau}}(t_j)\bigr)$ is eventually contained in
$\wbar{\eta_{\bar{\tau}}(H_0)}$.  Since $\eta_{\bar{\tau}}(H)$ is
dense in $H_{\bar{\tau}}$, this is sufficient to see that
$\wbar{\eta_{\bar{\tau}}(H_0)}$ is open.  \endpf

Immediately, we obtain a restriction theorem.

\begin{spinerest}\label{theo:spinerest}
If $H_0$ is an open subgroup of a locally compact group
$H$, then $\spineh|_{H_0}=\spineho$ and
$\rspineh|_{H_0}=\rspineho$.
\end{spinerest}

\proof It is clear that if $\tau\in\th$ then $\tau\cap H_0
=\{U\cap H_0:U\in\tau\}\in\tho$.
Thus if follows for each $\tau\in\th$ that $\falbh{\tau}|_{H_0}
\subseteq\falbho{{\tau\cap H_0}}$.  Hence it follows from
Theorem \ref{theo:directsum} that $\spineh|_{H_0}\subseteq\spineho$.
Now if $u\in\spineho$, from Theorem \ref{theo:directsum} we may find
for each $\tau\iin\tnqho$ a unique $u_\tau\iin\falbho{\tau}$ such that
\[
u=\sum_{\tau\in\tnqho}u_\tau\quad\wwhere\quad
\norm{u}=\sum_{\tau\in\tnqho}\norm{u_\tau}.
\]
Each $u_\tau$ corresponds to a function $\hat{u}_\tau\iin\fal{(H_0)_\tau}$.
By Lemma \ref{lem:exttop}, above, we can identify $(H_0)_\tau$
as an open subgroup of $H_{\bar{\tau}}$.
Let $\hat{v}_\tau\iin\fal{H_{\bar{\tau}}}$ be so 
$\hat{v}_\tau=\hat{u}_\tau$ on  $(H_0)_\tau$, and $\hat{v}_\tau=0$
off $(H_0)_\tau$.  Then let $v_\tau$ be the corresponding element of 
$\falbh{\bar{\tau}}$; we have that $\norm{v_\tau}=\norm{\hat{v}_\tau}=
\norm{u_\tau}$.  
If $v=\sum_{\tau\in\tnqho}v_\tau$, then $v\in\spineh$ with
$\norm{v}\leq\norm{u}$, and $v|_{H_0}=u$.  

We obtain the result for $\rspineh$ since both of the
operations $\tau\mapsto\tau\cap H_0$ on $\th$, and
$\tau\mapsto\bar{\tau}$ on $\tho$ create Hausdorff
topologies from Hausdorff topologies.  \endpf

We note that we cannot expect that $\spineh|_{H_0}=\spineho$
when $H_0$ is closed, but not open in $H$.  See (\ref{eq:restfail}).

Another immediate application of Lemma \ref{lem:exttop} is the following
theorem.  Let us recall that a {\it coset} $K$ of $G$ is a set
which satisfies the relation
\begin{equation}\label{eq:coset}
r,s,t\in C\qquad\text{implies}\qquad rs^{-1}t\in C.
\end{equation}
Moreover, $K^{-1}K$ and $KK^{-1}$ are then subgroups for which
$sK^{-1}K=KK^{-1}s=K$ for each $s\in K$. 
See \cite{ilies}, Proposition 1.1.  
Thus it follows that a non-empty intersection
of a family of cosets is itself a coset.

\begin{fsalidem}\label{theo:fsalidem}
Every idempotent in $\fsalh$ is an element of $\spineh$.
\end{fsalidem}

We note that for abelian groups this result is \cite[8.1.4]{taylor}.
There the result is proved independently of Cohen's idempotent
theorem \cite{cohen0}, and in fact is used to deduce it.
We use the idempotent theorem of Host \cite{host}.

\medskip
\proof Let us first consider an idempotent of the form $1_{H_0}$
where $H_0$ is an open subgroup of $H$.  If we let $\tauap^0$ denote
the topology on $H_0$ induced by its almost periodic compactification,
then $\tau_0=\wbar{\tauap^0}$, from Lemma \ref{lem:exttop},
is such that $1_{H_0}\in\falb{\tau_0}{H}$.  Note that for cosets we
have $1_{tH_0}=t\con 1_{H_0}$ and $1_{H_0t}=1_{H_0}\con t$ (where
$t\con u(s)=u(t^{-1}s)$ and $u\con t(s)=u(st^{-1})$ for $t,s\iin H$
and $u$ a function), so
$1_{tH_0},1_{H_0t}\in \falb{\tau_0}{H}$.

By \cite{host} each idempotent is of the form 
$1_Y$, where $Y\in\ocringh$, the ring generated by open cosets
of $H$.  For such $Y$, we can find open cosets $K_{ij}\subset L_i$,
for $j=1,\dots,n_i$ for which
$Y=\dot{\bigcup}_{i=1}^k\left(L_i\setdif\bigcup_{j=1}^{l_i}K_{ij}\right)$.
Let $H_1,\dots, H_n$ be a collection of open subgroups for which each of the 
non-empty cosets $L_i, K_{ij_1}\cap\dots\cap K_{ij_p}$, for $i=1,\dots,k$ and 
$1\leq j_1<\dots<j_p\leq l_i$, is a coset of one of these groups.
Letting $\tau_1,\dots,\tau_n$ be the topologies created as $\tau_0$ in the
above paragraph, we see that
\[
1_Y=\sum_{i=1}^k\left(1_{L_i}+
\sum_{p=1}^{l_i}(-1)^p\sum_{1\leq j_1<\dots<j_p\leq l_i} 
1_{K_{ij_1}\cap\dots\cap K_{ij_p}}\right)
\in\sum_{m=1}^n\falb{\tau_m}{H}
\]
and hence $1_Y\in\spineh$. \endpf

If $G$ and $H$ are both locally compact groups, then {\it continuous
piecewise affine} maps, $\alp:Y\subseteq H\to G$ where $Y\in\ocringh$,
were defined in \S \ref{ssec:homomorphisms}.  Given such a map
$\alp$, we define $\Phi_\alp:\falg\to\fsalh$ by
\[
\Phi_\alp u(s)=\begin{cases} u\bigl(\alp(s)\bigr) &\text{if }s\in Y \\
0 &\text{otherwise}. \end{cases}
\]
By \cite[Proposition 3.1]{ilies}, $\Phi_\alp$ is always a completely bounded
homomorphism.  These form the nicest class of homomorphisms from
$\falg$ to $\fsal{H}$, and were first studied for abelian groups
in \cite{cohen}.  

\begin{homorange}\label{theo:homorange}
If $G$ and $H$ are locally compact groups, $\alp:Y\subseteq H\to G$
is a continuous piecewise affine map and $\Phi_\alp:\falg\to\fsal{H}$
the homomorphism implemented by $\alp$, then there is a finite sequence
$\tau_1,\dots,\tau_n$ from $\tnqh$ such that
$\Phi_\alp\bigl(\falg\bigr)\subseteq\bigoplus_{i=1}^n\falbh{{\tau_i}}$.
Hence $\spineh$ contains the range of every such homomorphism.
\end{homorange}

\proof  First suppose that $Y$ is an open subgroup of $H$
and $\alp:Y\to G$ is a continuous homomorphism.  Then the topology
$\tau=\alp^{-1}(\tau_G)$ is a locally precompact group topology
on $Y$, whence by Lemma \ref{lem:exttop} there is a locally
precompact group topology $\bar{\tau}$ on $H$ which extends
$\tau$.  It then follows from Theorem \ref{theo:spinerest}
that $\falbh{ {\bar{\tau}} }|_Y=\falby{\tau}$.  Hence
$\Phi_\alp\bigl(\falg\bigr)=\falby{\tau}\subseteq\falbh{ {\bar{\tau}} }
\subseteq\falbh{{\bar{\tau}\vee\tauap}}$.

Now suppose that $Y$ is an open coset of $H$ and $\alp:Y\to G$ is a continuous
affine map.  Then if $H_0=Y^{-1}Y$ is an open subgroup of $H$ with 
$Y=s_0H_0$ for any $s_0\iin Y$.  Moreover
$\til{\alp}(s)=\alp(s_0)^{-1}\alp(s_0s)$ is a continuous
homomorphism on $H_0$.  Then for $u\iin\falg$ we see that
$\Phi_\alp u(s)=\alp(s_0)\con\Phi_{\til{\alp}}(s_0^{-1}\con u)$,
where $t\con v(s)=v(t^{-1}s)$ for a group element $t$ and function $v$.
Letting $\tau=\til{\alp}(\tau_G)$ on $Y$ and $\bar{\tau}$ it extension
to $H$, we have by translation invariance of $\falbh{{\bar{\tau}}}$
that $\Phi_\alp\bigl(\falg\bigr)=\Phi_{\til{\alp}}\bigl(\falg\bigr)
\subseteq\falbh{{\bar{\tau}}}\subseteq\falbh{{\bar{\tau}\vee\tauap}}$.

Now we suppose that $\alp:Y\to G$ is a general continuous piecewise affine
map.  Thus we can partition $Y=\dot{\bigcup}_{i=1}^mY_i$ where
each $Y_i\in\ocringh$ and there are open cosets $L_i\supseteq Y_i$
and continuous affine maps $\alp_i:L_i\to G$ so $\alp|_{Y_i}=\alp_i|_{Y_i}$
for $i=1,\dots,m$.  Thus we have that
\begin{equation}\label{eq:phialpdecomp}
\Phi_\alp=\sum_{i=1}^m\Phi_{\alp|_{Y_i}}.
\end{equation}
As in the proof of Theorem \ref{theo:fsalidem} above, we can find open
cosets $L_i$ and subsets $K_{i1},\dots,K_{in_i}$ of $L_i$
which are also open cosets such that
$Y_i=L_i\setdif\bigcup_{j=1}^{n_i}K_{ij}$.
But it then follows for each $i=1,\dots,m$ that
\begin{equation}\label{eq:phialpdecomp1}
\Phi_{\alp|_{Y_i}}=\Phi_{\alp_i}+
\sum_{p=1}^{n_i}(-1)^p\sum_{1\leq j_1<\dots <j_p\leq n_i}
\Phi_{\alp_i |_{K_{ij_1}\cap\dots\cap K_{ij_p}}}
\end{equation}
where we let $\Phi_{\alp_i |_{K_{ij_1}\cap\dots\cap K_{ij_p}}}=0$ if
$K_{ij_1}\cap\dots\cap K_{ij_p}=\varnothing$.
Since each $\alp_i$, $\alp_i|_{K_{ij_1}\cap\dots\cap K_{ij_p}}$ 
is affine, we see that there are topologies 
$\tau_i,\tau_{i;j_1},\dots,\tau_{i;1,\dots,n_i}$ in $\tnqh$ such that
\[
\Phi_{\alp_i}\bigl(\falg\bigr)\subseteq\falbh{\tau_i}\quad\aand\quad
\Phi_{\alp_i|_{K_{ij_1}\cap\dots\cap K_{ij_p}}}
\bigl(\falg\bigr)\subseteq\falbh{\tau_{i;j_1,\dots,j_p}}
\]
for $1\leq j_1<\dots j_p\leq n_i$ and $p=1,\dots,n_i$.  
Hence combining (\ref{eq:phialpdecomp}) and
(\ref{eq:phialpdecomp1}), we obtain
\[
\Phi_\alp\bigl(\falg\bigr)\subset
\sum_{i=1}^m\left(\falbh{\tau_i}+\sum_{j=1}^{n_i}
\sum_{1\leq j_1<\dots <j_p\leq n_i}\falbh{\tau_{i;j_1,\dots,j_p}}\right).
\]
Let $\tau_1,\dots,\tau_n$ be an enumeration of
of this set of topologies and we are done. \endpf

\begin{homorange1}\label{theo:homorange1}
Let $G$ and $H$ be locally compact groups, such that $G$ is amenable.

{\bf (i)} If $\Phi:\falg\to\fsal{H}$ is a completely bounded homomorphism
then there is a finite sequence
$\tau_1,\dots,\tau_n$ from $\tnqh$ such that
\[
\Phi\bigl(\falg\bigr)\subseteq\bigoplus_{i=1}^n\falbh{{\tau_i}}.
\]

{\bf(ii)} If $\Phi:\spineg\to\fsal{H}$ is a completely bounded homomorphism
then $\Phi\bigl(\spineg\bigr)\subseteq\spineh$.
\end{homorange1}

\proof (i) This follows from the theorem above and the main result of
\cite{ilies}.  

(ii)  By Theorem \ref{theo:directsum} we have
\[
\spineg=\ell^1\text{-}\!\!\!\!\bigoplus_{\tau\in\tnqg}\falbg{\tau}\cong
\ell^1\text{-}\!\!\!\!\bigoplus_{\tau\in\tnqg}\fal{G_\tau}.
\]
Each $G_\tau$ is amenable by \cite[Proposition 1.2.1]{rundeB},
so $\Phi|_{\falbg{\tau}}$ thus induces a completely bounded homomorphism
$\Phi_\tau:\fal{G_\tau}\to\fsal{H}$.  By (i) $\Phi_\tau\bigl(\fal{G_\tau}\bigr)
\subseteq\spineh$, so $\Phi\bigl(\falbg{\tau}\bigr)\subseteq\spineh$.
The result follows.  \endpf

Note that for an abelian $H$, it follows from \cite[8.2.4]{taylor},
and Theorem \ref{theo:abdual}, below,
that every invertible element $u\iin\fsalh$ is of the form
$u=v_1\dots v_ne^w$ where $v_i\in\falb{\tau_i}{H}\oplus\Cee 1$ 
for some $\tau_i\iin\tnqh$.  It would be interesting to determine
if this is true in general.  It seems likely that completing the
investigation begun in \cite{walter} might resolve this question.

\section{Abelian Groups}\label{sec:abgroup}

In this section we will let $G$ be a locally compact {\it abelian} group.  
We let
$\what{G}$ denote the {\it dual group} and $\tau_{\what{G}}$ its topology.  
If $\tau\in\tg$ then we let
$\what{G_\tau}$ denote the dual group of $G_\tau$ and $\hat{\tau}$
the topology on $\what{G_\tau}$.

Non-quotient topologies on abelian groups admit a particularly
appealing description.  We note that abelian groups are maximally
almost periodic since $\what{G}$ separates points on $G$.  Moreover,
we have that 
\begin{equation}\label{eq:abelapc}
\mapg=\what{(\what{G}_d)}
\end{equation}
where $\what{G}_d$ is the
topological group $\what{G}$ made discrete.

\begin{abdual}\label{theo:abdual}
Let $\tau\in\tg$.  Then $\tau\in\tnqg$ if and only if  
$\what{G_\tau}=\what{G}$ and $\hat{\tau}\supseteq\tau_{\what{G}}$.
Moreover, any locally compact group topology
on $\what{G}$ which is finer than
$\tau_{\what{G}}$ is of the form $\hat{\tau}$ for some $\tau\iin\tnqg$.
\end{abdual}

\proof First, we see that if $\tau\in\tg$ then the homomorphism
$\eta_\tau:G\to G_\tau$ has dense range, so the dual map
$\hat{\eta}_\tau:\what{G_\tau}\to\what{G}$, given by
$\hat{\eta}_\tau(\sig)=\sig\comp\eta_\tau$, is injective.

Now suppose that $\tau\in\tnqg$ so $\tau=\tau\vee\tauap$ by Theorem
\ref{theo:nonquot} (iii).  Thus we get the diagram
\begin{equation}\label{eq:abelnq}
\begin{CD}
G@>{\eta_{\tau\vee\tauap}}>>G_{\tau\vee\tauap}@>{\iota}>>
G_\tau\cross \mapg
\end{CD}
\end{equation}
where $\iota$ is the injection.  Then (\ref{eq:abelapc}) and the Pontryagin 
duality
theorem tell us that $\what{\mapg}\cong\what{G}_d$.  Hence (\ref{eq:abelnq})
admits dual diagram
\[
\begin{CD}
\what{G_\tau}\cross\what{G}_d@>{(\sig_\tau,\sig)\mapsto\sig_\tau\cross\sig}>>
\what{G_\tau\cross\mapg}@>{r}>>\what{G_{\tau\vee\tauap}}
@>{\hat{\eta}_{\tau\vee\tauap}}>>\what{G}
\end{CD}
\]
where $r$ is the restriction map, which is surjective
by \cite[24.12]{hewittrI}, and $\sig_\tau\cross\sig$ is the 
Kroenecker product of the characters $\sig_\tau$ and $\sig$.
If $\sig\in\what{G}$ then for any $s\iin G$ we have that
$1_\tau\cross\sig\bigl(\eta_{\tau\vee\tauap}(s)\bigr)=\sig(s)$
which shows that $\hat{\eta}_{\tau\vee\tauap}$ is surjective.

Thus we see that $\hat{\eta}_\tau:\what{G}\to\what{G_\tau}$ is continuous and 
bijective.  It follows that $\hat{\tau}\supseteq\tau_{\what{G}}$.

Now let us suppose that $\what{G_\tau}=\what{G}$.  Then
there is a natural bicontinuous bijection between the diagonal
subgroup of $\what{G_\tau}\cross\what{G}_d$ and $\what{G_\tau}$.
However, the diagonal subgroup of $\what{G_\tau}\cross\what{G}_d$
is exactly $\what{G_{\tau\vee\tauap}}$.  Thus
$\what{G_\tau}\cong\what{G_{\tau\vee\tauap}}$, bicontinuously.  Using
Pontryagin duality we see that $G_\tau\cong G_{\tau\vee\tauap}$,
bicontinuously, so $\tau=\tau\vee\tauap$, whence
$\tau\in\tnqg$ by Theorem \ref{theo:nonquot} (iii). 

Finally, if $\til{\tau}$ is any locally compact group topology on $\what{G}$
which is finer than $\tau_{\what{G}}$, then the continuous bijection
$\id_{\what{G}}:(\what{G},\til{\tau})\to \what{G}$ has dual map 
$\eta:G\to\what{(\what{G},\til{\tau})}$, and induces a 
locally precompact topology
$\tau$ on $G$ for which $\til{\tau}=\hat{\tau}$.  It follows
from above that $\tau\in\tnqg$.  \endpf

It is worth noting the dual to our lattice operation on $\tnqg$.

\begin{abdual1}\label{cor:abdual1}
If $\tau_1,\tau_2\in\tnqg$ then 
$\what{\tau_1\vee\tau_2}=\hat{\tau}_1\cap\hat{\tau}_2$
on $\what{G}$.
\end{abdual1}

\proof For $j=1,2$ the homomorphism 
$\eta^{\tau_1\vee\tau_2}_{\tau_j}:G_{\tau_1\vee\tau_2}
\to G_{\tau_j}$ has continuous surjective dual map
$\hat{\eta}^{\tau_1\vee\tau_2}_{\tau_j}:(\what{G},\hat{\tau}_j)\to
(\what{G},\what{\tau_1\vee\tau_2})$, so 
$\what{\tau_1\vee\tau_2}\subseteq\tau_j$.
Hence $\what{\tau_1\vee\tau_2}\subseteq\hat{\tau}_1\cap\hat{\tau}_2$.
However, $\hat{\tau}_1\cap\hat{\tau}_2$ is evidently a locally precompact
topology on $\what{G}$ whose dual topology 
$\what{\hat{\tau}_1\cap\hat{\tau}_2}$,
which we may consider as a topology on $G$,
dominates each of $\tau_1,\tau_2$, and is dominated by $\tau_1\vee\tau_2$.  
Hence $\what{\hat{\tau}_1\cap\hat{\tau}_2}=\tau_1\vee\tau_2$
and by Pontryagin duality we must have 
$\what{\tau_1\vee\tau_2}=\hat{\tau}_1\cap\hat{\tau}_2$. \endpf

\section{Examples}

We close this paper with some examples to indicate the scope
of $\tnqg$ and hence $\spineg$ and its spectrum $G^*$.

First note that if $G$ is compact, then $\tnqg=\{\tau_G\}$,
since every coarser group topology is clearly a quotient of $\tau_G$.

We will compute several slightly less trivial examples below, where
$\tnqg=\{\tauap,\tau_G\}$.  As trivial as this seems, it still give rise
to an interesting spine compactification of $G$:  
\[
G^*=\mapg\sqcup G
\]
by (\ref{eq:tnqgcompletespec}).  As in Theorem \ref{theo:spectopology}
the topology is given by letting
$G$ be an open subgroup of $G^*$, and for any point in
$\mapg$, a basic open neighbourhood is of the form
$U\sqcup W$, where $U$ is open in $\mapg$ and
$G\setdif W$ is compact in $G$.  The semigroup multiplication
is given by allowing $\mapg$ and $G$ to each be subgroups
while for $s\iin\mapg$ and $t\iin G$ we let 
$st=s\etaap(t)$ and $ts=\etaap(t)s$.  Hence 
$\mapg$ is a closed ideal in $G^*$.

We explicitly compute $\tnqg$ for vector groups
and certain semi-direct product Lie groups in
\S \ref{ssec:liegroups}.

In the last two subsections we give some hints at the scope and
complexity of $\tnqg$ for algebraic groups and free groups.

Since it will prove useful below, let us 
note the following fact.

\begin{quotgroup}\label{lem:quotgroup}
If $G$ is a locally compact group, $K$ a compact normal subgroup
with quotient map $q_K:G\to G/K$, and $\tau\in\tnqg$,
then $\bar{\tau}=\{U\subseteq G/K:q_K^{-1}(U)\in\tau\}\in\tnq{G/K}$.
\end{quotgroup}

\proof  First note that $\eta_\tau(K)$ is a compact normal subgroup
of $G_\tau$.  Thus $\eta_\tau$ induces a continuous homomorphism
$\bar{\eta}_\tau:G/K\to G_\tau/\eta_\tau(K)$, and
$\bar{\tau}=\bar{\eta}_\tau^{-1}(\tau_{G_\tau/\eta_\tau(K)})$
and is thus locally precompact.  Let $(\mapg,\etaap)$ denote
the almost periodic compactification of $G$.  Then, by similar
reasoning as above, we have that $\etaap(K)$ is bicontinuously
isomorphic to $K$ and $\mapg/\etaap(K)\cong(G/K)^{ap}$.
Thus if $U\in\tauap(G/K)$ -- the topology on $G/K$ induced by its
almost periodic compactification -- there is $V\iin\tauap$ such that
$q_K^{-1}(U)=VK$, so $q_K^{-1}(U)\in\tauap\subseteq\tau$.
Hence $\tauap(G/K)\subseteq\bar{\tau}$ so $\bar{\tau}\in\tnq{G/K}$.
\endpf

\subsection{Some Lie Groups}  \label{ssec:liegroups}
We begin with the most basic class of
examples.  We will use additive notation when dealing with the 
integers $\Zee$,
real numbers $\Ree$, and for any vector group.

\begin{reeandtee}\label{prop:reeandtee}
If $G$ admits a compact normal subgroup $K$ for which
$G/K$ is bicontinuously isomorphic to either
$\Ree$ or $\Zee$,  then $\tnqg=\{\tauap,\tau_G\}$.
\end{reeandtee}

We note that if $K$ is trivial, then Theorem \ref{theo:abdual}
tells us that on $\Ree\cong\what{\Ree}$, or $\Tee\cong\what{\Zee}$,
the only locally compact group topology strictly finer than
the usual topology is the discrete topology.  This fact
seems well-known.  See, for example, \cite[Section 2]{rickert}.

\medskip
\proof If $\tau\in\tnqg$, then by Lemma \ref{lem:quotgroup}
we obtain a non-quotient topology $\bar{\tau}$
on $G/K$.  Then by \cite[9.1]{hewittrI}, the continuous homomorphism
$\eta_{\bar{\tau}}:G/K\to (G/K)_{\bar{\tau}}$ satisfies either that
$\eta_{\bar{\tau}}(G/K)$ is bicontinuously isomorphic to
$G/K$, or has compact closure.  In the first case,
$\tau=\tau_G$; in the second, $\tau=\tauap$. \endpf

We note that the discrete reals $\Ree_d$ will present a much more
complicated structure.  We will see in the next section that the additive
rationals $\Que$, which is an open subgroup, presents a non-trivial
structure.  

We generalise Proposition \ref{prop:reeandtee} to vector groups.
The following is suggested by \cite[Lemma 2.6]{rickert}.
We recall that $\what{\Ree^n}\cong\Ree^n$ for any non-negative integer
$n$.   

\begin{vectorgroup}\label{theo:vectorgroup}
If $n\geq 2$, let $\fL_n$ denote the subspace lattice of $\Ree^n$.
For each $L\in\fL_n$ we let $\hat{\tau}_L$ be the topology on $\Ree^n$
formed by letting $L$ be an open subgroup and $\hat{\tau}_L\cap L$
be the usual topology.  Then $\hat{\tau}_L$ is locally compact on 
$\Ree^n$ and finer than $\tau_{\Ree^n}$.  Moreover if we let
$\tau_L$ be the locally precompact non-quotient
topology on $\Ree^n$ for which $\what{\tau_L}=\hat{\tau}_L$ 
(in the notation of Section \ref{sec:abgroup}), then
\[
\tnq{\Ree^n}=\{\tau_L:L\in\fL_n\}
\]
with semilattice structure
\[
\tau_{L_1}\vee\tau_{L_2}=\tau_{L_1+ L_2}.
\]
Hence $\tnq{\Ree^n}$ is isomorphic to the upper semilattice $\fL_n$.
\end{vectorgroup}

Note that in the notation above $\tau_{\Ree^n}$ is the usual topology
while $\tau_{\{0\}}=\tauap$.  The above result extends very easily
to all groups $G$ for which $G/K\cong\Ree^n$ for a compact normal
subgroup $K$.  Note that in $\fL_n$, every hereditary directed subset
in $\tnqg$ is principal.  Hence by Proposition \ref{prop:hdprincipal} 
is a lattice with 
\begin{equation}\label{eq:lsubspacelattice}
\tau_{L_1}\wedge\tau_{L_2}=\tau_{L_1\cap L_2}.
\end{equation}
Indeed, it is clear that $\tau_{L_1\cap L_2}\subseteq\tau_{L_j}$ for $j=1,2$.
If $L_0\iin\fL_n$ satisfies $\tau_{L_0}\subset\tau_{L_j}$ for each $j=1,2$, 
then $L_0\subseteq L_j$ for each $j$, so $L_0\subseteq L_1\cap L_2$.

\medskip
\proof  It is plain to see that if $L\in\fL_n$, then
$\hat{\tau}_L$ forms a locally compact group topology on $\Ree^n$.
Moreover, if $L'$ is a linear complement,
then since $L$ is open, $\{0\}=L\cap L'\in\hat{\tau}_L\cap L'$;
so $\hat{\tau}_L\cap L'$ is the discrete topology.  It follows
from the direct sum decomposition $\Ree^n=L+L'$ that $(s,t)\mapsto
s+t$ is a bicontinuous isomorphism from $L\cross L'_d$ to 
$(\Ree^n,\til{\tau}_L)$.  It is obvious
that $\hat{\tau}_L$ is finer that $\tau_{\Ree^n}$.  It thus follows from
Theorem \ref{theo:abdual} that there is $\tau_L\iin\tnq{\Ree^n}$ for
which $\hat{\tau}_L=\what{\tau_L}$.

Let us now check the lattice operations.  If $L_1,L_2\in\fL_n$
then it is clear that $\hat{\tau}_{L_j}\supseteq\hat{\tau}_{L_1+L_2}$
for $j=1,2$.  On the other hand, if $U\in
\hat{\tau}_{L_1}\cap\hat{\tau}_{L_2}$ and $s\in U$ then
find a neighbourhood of the identity $V$ with 
$s+2V\subseteq U$.  Then
$W=V\cap L_1+V\cap L_2\in\hat{\tau}_{L_1+L_2}$, since $L_1+L_2$
is a linear quotient of $L_1\cross L_2$,
and $s+W\subseteq s+2V\subseteq U$. 
Hence $U\in\hat{\tau}_{L_1+L_2}$.
Thus $\hat{\tau}_{L_1}\cap\hat{\tau}_{L_2}=
\hat{\tau}_{L_1+L_2}$ and it follows from Corollary \ref{cor:abdual1}
that $\tau_{L_1}\vee\tau_{L_2}=\tau_{L_1+L_2}$.

Now let $\hat{\tau}$ be any locally compact group topology on $\Ree^n$
which is finer than $\tau_{\Ree^n}$.  It remains to show that $\hat{\tau}=
\hat{\tau}_L$ for some $L\iin\fL_n$.  Let $U$ be a symmetric
open neighbourhood of $0$ in $(\Ree^n,\hat{\tau})$ with compact
closure.  Then $H=\bigcup_{n=-\infty}^\infty nU$ is a compactly
generated open subgroup of $(\Ree^n,\hat{\tau})$.  It follows
\cite[9.8]{hewittrI} that there are non-negative integers
$l$ and $k$ and a compact abelian group $K$ such that
$H\cong\Ree^l\cross\Zee^k\cross K$ as topological groups.
Since $\id:(\Ree^n,\hat{\tau})\to\Ree^n$ is continuous,
the injection $H\hookrightarrow(\Ree^n,\hat{\tau})$ induces
a continuous homomorphism $\iota:\Ree^l\cross\Zee^k\cross K
\to\Ree^n$.  It follows that $K$ is trivial and $0\leq l\leq n$.
Let $L=\iota(\Ree^l)$.  Then $L$ is an open subgroup of $H$, and
hence an open subgroup of $(\Ree^n,\hat{\tau})$.  Moreover,
$\hat{\tau}\cap L=\tau_{\Ree^l}$.  This implies that $\hat{\tau}
=\hat{\tau}_L$. \endpf

We can compute that for each $L\iin \fL_n$, where 
$L'$ is a complementary subspace, that
$\Ree^n_{\tau_L}\cong L\cross (L')^{ap}$.
Hence by (\ref{eq:tnqgcompletespec}) we have
\[
(\Ree^n)^*\cong\bigsqcup_{L\in\fL_n} L\cross (L')^{ap}.
\]
Moreover, we can use (\ref{eq:gstarmult1})
to obtain the semigroup operation.  If $s=(s_L,s_{ap})\iin L\cross(L')^{ap}$
and $t=(t_M,t_{ap})\iin M\cross(M')^{ap}$, we let $p:\Ree^n\to L\cap M$
be the projection relative to the decomposition $\Ree^n:L\cap M+(L'+M')$
and $p'=\id-p$.  Then the semigroup product of $s$ and $t$ is given
as in (\ref{eq:gstarmult1})
\[
\biggl(p(s_L)+p(t_M),\bigl(\etaap^L\comp p'(s_L),s_{ap}\bigr)+
\bigl(\etaap^M\comp p'(t_M),t_{ap}\bigr)\biggr)
\in L\cap M\cross (L'+M')^{ap}
\]
where $\etaap^L:p'(L)\to p'(L)^{ap}$ is the almost periodic compactification
map, and $\etaap^M$ is defined similarly.  The topology can be described
by Theorem \ref{theo:spectopology}.  To observe some of its features
let for each $k=1,\dots,n$ 
\[
\fL_{n,k}=\{L\in\fL_n:\dim L=k\}\quad\aand\quad
(\Ree^n)^*_k=\bigsqcup_{L\in\fL_{n,k}}L\cross(L')^{ap}.
\]
Then each set $\bigsqcup_{k=1}^l(\Ree^n)^*_k$ is a closed ideal.
Also, in each of the subgroupoids
$(\Ree^n)^*_k$ (a Clifford semigroup is a {\it groupoid} if we forget to 
multiply between component subgroups; see 
\cite[Proposition 1.0.1]{pattersonG}), in the relativised topology, 
each of the components $L\cross (L')^{ap}$ is an open subgroup.

We can also compute the lattice of non-quotient topologies for
the finitely generated free abelian group $\Zee^n$.  Recall
that $\what{\Zee^n}\cong\Tee^n$.

\begin{lattice}\label{cor:lattice}
The semilattice $\tnq{\Zee^n}$ is isomorphic to $\fL_n$.
\end{lattice}

\proof As in the proof of Theorem \ref{theo:vectorgroup} we
will identify locally compact group topologies on $\Tee^n\cong\what{\Zee^n}$
which are finer than $\tau_{\Tee^n}$ and appeal to Theorem \ref{theo:abdual}.
We refer the reader to \cite{warner} for pertinent definitions
and elementary results in Lie groups.

Let $q:\Ree\to\Tee$ be the quotient map given by $q(t)=e^{i2\pi t}$.
Then the $n$-fold Kroenecker product map $q^n:\Ree^n\to\Tee^n$ is
a smooth quotient map -- the universal covering map.  
Thus for distinct $L,M\iin\fL_n$ we have that
$q^n(L)$ and $q^n(M)$ are distinct Lie subgroups of $\Tee^n$.
Then $\hat{\tau}_{\{0\}}$ is the discrete topology.
For a given $L\iin\fL_n\setdif\{0\}$ 
we have that there is a Lie group isomorphism
$q^n(L)\cong L/(L\cap\ker q^n)\cong\Ree^l\cross\Tee^m$, 
where either $0<l$, $0\leq m$ and
$l+m<n$, or $l=0$ and $0<m\leq n$.  We then let 
$\hat{\tau}_L$ be the coarsest topology on $\Tee^n$ which makes
$q^n(L)$ an open subgroup which is bicontinuously isomorphic to
$\Ree^l\cross\Tee^m$.  As in the proof of the theorem above
we have that $\hat{\tau}_L\cap\hat{\tau}_M=\hat{\tau}_{L+M}$.
Thus if $\tau_L$ is the topology on $\Zee^n$ for which $\what{\tau_L}
=\hat{\tau}_L$, then $\{\tau_L:L\in\fL_n\}$ forms a subsemilattice
of $\tnq{\Zee^n}$.

It remains to show that every locally compact group topology $\hat{\tau}$ on
$\Tee^n$, finer than $\tau_{\Tee^n}$, is of the form $\hat{\tau}_L$
for some $L\iin\fL_n$.  As in the proof of the above theorem we
find a compactly generated open subgroup $H$.  Then $H\cong
\Zee^k\cross\Ree^l\cross K$ where $K$ is a compact abelian group.
Since $K$ must be a compact, hence closed, subgroup of $\Tee^n$,
it is necessarily a Lie subgroup and thus has an open subgroup of
the form $\Tee^m$, where $m\leq n$.  The open subgroup of $H$
corresponding to $\Ree^l\cross\Tee^m$ injects continuously in to 
$\Tee^n$, and hence must be one of the Lie subgroups of the form
$q^n(L)$, where $L\in\fL_n$, as above.  \endpf

Let us illustrate the above result by computing
$\Zee^2_\tau$ for each $\tau\iin\tnq{\Zee^2}$.  
We have that $\tnq{\Zee^2}\cong\fL_2=\bigl\{\Ree^2,L_\theta,\{0\}:
0\leq\theta<\pi\bigr\}$ where for each $\theta$, 
$L_\theta$ is the subspace generated by $(\cos\theta,
\sin\theta)$.  If $\tan\theta$ is rational (or $\theta=\pi/2$) then
$q^2(L_\theta)\cong\Tee$; while otherwise 
$q^2(L_\theta)\cong\Ree$,  a skew-line in $\Tee^2$.
Letting $\tau_\theta=\tau_{L_\theta}$ we have that
\[
\Zee^2_{\tau_\theta}=\wbar{\left\{\bigl(n\cos\theta+m\sin\theta,\eta_{ap}(n,m)
\bigr):(n,m)\in\Zee^2\right\}}\subset\Ree\cross(\Zee^2)^{ap}.
\]
Indeed, if $\eta_\theta:\Zee^2\to\Ree$ is given by $\eta_\theta(n,m)
=nc+ms$ ($c=\cos\theta,s=\sin\theta$), $\xi\in\Ree$ and $\hat{\xi}:\Ree\to\Tee$
is the character $\hat{\xi}(t)=e^{i2\pi\xi t}$, then
$\hat{\xi}\comp\eta_\theta(n,m)=e^{i2\pi\xi(nc+ms)}=q(\xi c)^nq(\xi s)^m$.
Thus $\xi\mapsto q^2(\xi c,\xi s)$ is the dual map to $\eta_\theta$, and has its
range in $q^2(L_\theta)$.
We then can see that $\tau_\theta=\eta_\theta^{-1}(\tau_\Ree)\vee\tauap$.
Note that if $\tan\theta=s/c$ is irrational, then $\eta_\theta(\Zee^2)$ is 
dense in $\Ree$; otherwise it is a closed subgroup.

\medskip
Let us now embark on computing $\tnqg$ for any group $G$ for which
the following formula is satisfied
\begin{equation}\label{eq:mwapdec}
\fsalg=\falbg{\tauap}\oplus^{\ell^1}\fsalog
\end{equation}
where $\fsalog$ is the closed subalgebra of $\fsalg$ of functions
vanishing at infinity.
This class of groups includes all {\it minimally weakly almost periodic}
groups, i.e., those for which the algebra of weakly almost periodic functions 
$\mathcal{WAP}(G)$ is as small as it can be, namely
\[
\mathcal{WAP}(G)=\mathcal{AP}(G)\oplus\mathcal{C}_0(G).
\]
See \cite{chou} for more on this.
The class of minimally almost periodic Lie groups includes 
$\mathrm{SL}_2(\Ree)$
and the motion groups $\mathrm{M}(n)=\mathrm{SO}(n)\!\ltimes\!\Ree^n$.

These types of groups admit small spines.

\begin{mwap}\label{theo:mwap}
If $G$ satisfies (\ref{eq:mwapdec}), then
$\tnqg=\{\tauap,\tau_G\}$.  Hence $\spineg=\falbg{\tauap}\oplus^{\ell^1}\falg$.
\end{mwap}

\proof If $\tau\in\tnqg$ and $\fA=\faltg\cap\fsalog\not=\{0\}$, then
$\fA$ is isomorphic to a translation invariant $*$-closed subalgebra
of $\faltg$.  Hence by Lemma \ref{lem:subalgebraq}, there is
a quotient topology $\tau_0$ of $\tau$ for which $\fA=\falbg{\tau_0}$.
Then $\eta_{\tau_0}:G\to G_{\tau_0}$ must be a proper map.  Indeed,
if $K\subset G_{\tau_0}$ and is compact, find $\hat{u}\iin\fal{G_{\tau_0}}$
so that $\hat{u}|_K=1$.  Then $\eta_{\tau_0}^{-1}(K)$ is a closed
subset of $G$, and if $u=\hat{u}\comp\eta_{\tau_0}$ then
$u|_{\eta_{\tau_0}^{-1}(K)}=1$ and $u\in\fsalog$, so
$\eta_{\tau_0}^{-1}(K)$ is compact.  In other words, $\eta_{\tau_0}$
is proper, so $\tau_0$ is a quotient of $\tau_G$.  
Hence by Theorem \ref{theo:nonquot} (ii), $\tau=\tau_G$.

Thus if $\tau\in\tnqg$ and $\tau\not=\tau_G$, then by
(\ref{eq:mwapdec}) $\faltg\subseteq\falbg{\tauap}$,
which forces $\tau=\tauap$.  \endpf

We note that the above result can also be derived for
$\mathrm{SL}_2(\Ree)$ by noting that it
is a {\it totally minimal group} (see \cite[7.4.1]{dikranjanps})
with only one normal subgroup $Z=\{-I,I\}$, which is compact.
Moreover $\mathrm{SL}_2(\Ree)^{ap}=1$.  It follows that
the spine compactification coincides with the one-point
compactification.

From the proof of Theorem
\ref{theo:mwap} we see that $\spineg\cap\fsalog=\falg$.  Thus it follows
from the fact that $\fsalog\not=\falg$ for many unimodular groups,
by \cite{figatalamanca}, that $\spineg\not=\fsalg$ for such groups.
It is noted in \cite{inoue1} that $\spineg\not=\fsalg$ for any non-compact
abelian group.

In \cite{walter76} and in \cite[Corollary 2.1]{mauceri}, 
it is shown that the ``$p$-adic motion group'' 
$G=\Tee_p\!\ltimes\!\Que_p$, where $\Que_p$ is the $p$-adic numbers ($p$ prime)
and $\Tee_p=\{r\in\Que_p:|r|_p=1\}$ (see Section \ref{ssec:alggroups}, below),
satisfies a strong form of (\ref{eq:mwapdec}):
$\fsalg=\falbg{\tauap}\oplus^{\ell^1}\falg$.
Thus $\spineg=\fsalg$ for this non-compact group.

\medskip
Let us close this section by examining the $ax+b$-group
$H=\{(a,b):a,b\in\Ree,a>0\}$ with multiplication
$(a,b)(a_1,b_1)=(aa_1,ab_1+b)$.  If $j:H\to \Ree$ is the homomorphism
given by $j(a,b)=\log a$, then $j$ is continuous with 
$\ker j=\{(1,b):b\in\Ree\}\cong\Ree$.  In \cite{gelfandn}
the family of irreducible representations is computed;
each is a character whose kernel includes $\ker j$.  
It thus follows that
that the maximal almost periodic compactification map 
$\etaap$ factors through $j$, and hence $H^{ap}\cong\Ree^{ap}$
Moreover, in \cite[Th\'{e}or\`{e}me 7]{khalil} it is shown that
\[
\fsalh=\bigl(\fsal{\Ree}\comp j\bigr)\oplus^{\ell^1}\fal{H}.
\]

\begin{axpb}\label{theo:axpb}
If $H$ is the $ax+b$-group, as above, then $\tnqh
=\{\tauap,\til{\tau}_\Ree,\tau_H\}$, where $\til{\tau}_\Ree=j^{-1}(\tau_\Ree)$.
Hence
\[
\spineh=
\falbg{\tauap}\oplus^{\ell^1}\falb{\til{\tau}_\Ree}{H}\oplus^{\ell^1}\fal{H}
\cong\fal{\Ree^{ap}}\oplus^{\ell^1}\fal{\Ree}\oplus^{\ell^1}\fal{H}.
\]
\end{axpb}

\proof In a manner 
similar to the proof of Theorem \ref{theo:mwap}, we can show that
if $\tau\in\tnqh$ and $\tau\not=\tau_H$, then $\falb{\tau}{H}\subset
\fsal{\Ree}\comp j$.  Hence by Proposition \ref{prop:containment}
we have that $\tau\subseteq\til{\tau}_\Ree$.  Thus 
$\eta^{\til{\tau}_\Ree}_{\tau}:H_{\til{\tau}_\Ree}\to H_\tau$
induces a non-quotient topology on $\Ree\cong H_{\til{\tau}_\Ree}$.
It then follows Proposition \ref{prop:reeandtee} that
$\tau$ is either $\til{\tau}_\Ree$ or $\tauap$.  \endpf

We thus deduce that 
\begin{equation}\label{eq:restfail}
\spineh|_{\ker j}\cong\Cee 1\oplus^{\ell^1}\fal{\Ree}\not\cong\spine{\Ree}.
\end{equation}
Also, $\rspine{H}=\falh$, so 
\[
\rspineh|_{\ker j}\cong\fal{\Ree}\not\cong\rspine{\Ree}.
\]

We note that $\tnqhh=\{\tau_H\}$ and $\rspineh=\fal{H}$. Moreover
we get spine compactification $H^*=\Ree^{ap}\sqcup\Ree\sqcup H$.
The topology admits $H$ and $\Ree\sqcup H$ as open subsets.
For any $r\iin \Ree^{ap}$ a basic open neighbourhood is of the form
$U\sqcup\bigl((-\infty,-a)\cup(a,\infty)\bigr)\sqcup W$ where 
$U$ is an open neighbourhood of $r$ in $\Ree^{ap}$, $a\geq 0$ and
$H\setdif W$ is compact; for $s\in\Ree$ a basic open neighbourhood
is $(s-\eps,s+\eps)\sqcup W$ where $\eps>0$ and $W$ is as above.
If $r\in\Ree^{ap}$, $s\in\Ree$ and $t\in H$ then we have multiplications
\[
rs=r\etaap(s),\quad st=sj(t)\;\aand\;rt=r\etaap\comp j(t)
\]
with similar formulas for $sr$, $ts$ and $tr$, where $\etaap:\Ree\to\Ree^{ap}$
is the compactification map.

\subsection{Some Algebraic Groups}\label{ssec:alggroups}

Let us first look at some topologies on the discrete additive
rational numbers $\Que$.
We let $\tau_\Que$ denote the discrete topology.  We
let $\til{\tau}_\Ree$ denote the order topology.  Finally, we 
let $\Pee$ denote the set of primes in $\En$ and 
for any $p\iin\Pee$, let $O_p=\Zee[1/q:q\in\Pee,q\not=p]$.
Then let $\til{\tau}_{\Que_p}$ denote the 
group topology for which $\{p^kO_p:k\in\Zee\}$
is a base at $0$.  Then $\Que_{\til{\tau}_\Ree}\cong\Ree$;
and $\Que_{\til{\tau}_{\Que_p}}\cong\Que_p$, the {\it $p$-adic
numbers}.  (Note that one of our main references \cite{hewittrI}
uses the notation $\Omega_p=\Que_p$.)  It is well-known that $\Que=\Que_d$,
$\Ree$ and $\Que_p$ ($p\in\Pee$) represent all of the locally
compact fields in which $\Que$ is a dense subfield;
see \cite[Theorem 4-12]{ramakrishnanv}, for example.  In fact,
$\Ree$ is the metrical completion of $\Que$ via the metric
$\rho_\Ree(r,s)=|r-s|$, where $|\cdot|$ is the usual order-induced absolute
value; and $\Que_p$ is the metrical completion of $\Que$ via
the metric $\rho_p(r,s)=|r-s|_p$, where $|\cdot|_p$ is the
{\it $p$-adic norm} which is given by
\[
|r|_p=p^{-\nu_p(r)}\;\wwhere\; 
\nu_p(r)=\sup\left\{\nu\in\Zee:r\in p^\nu O_p\right\}.
\]
Let us identify $\Que$ with its natural copies in $\Ree$
and each $\Que_p$.  It follows from the weak approximation theorem
from number theory (see \cite[Section 6]{cassels}, for eaxmple),
that for any collection $p_1,\dots,p_n$ of distinct primes we have that
the ``diagonal" subgroup $\{(r,\dots,r):r\in\Que\}\subset\Que^{n+1}$ is
dense in $\Ree\cross\Que_{p_1}\cross\dots\cross\Que_{p_n}$.
In other words
\[
\Que_{\til{\tau}_\Ree\vee\til{\tau}_{\Que_{p_1}}\vee\dots\vee
\til{\tau}_{\Que_{p_n}}}
\cong\Ree\cross\Que_{p_1}\cross\dots\cross\Que_{p_n}.
\]
(The second author is grateful to Eugene Eisenstein for pointing this
fact out to him.)  Let us also note the following, in a analogy with
Proposition \ref{prop:reeandtee}.

\begin{padic}\label{prop:padic}
For any $p\iin\Pee$, $\tnq{\Que_p}=\{\tauap,\tau_{\Que_p}\}$
\end{padic}

\proof We recall the well-known fact that $\Zee_p=\{r\in\Que_p:|r|_p\leq 1\}
=\wbar{O_p}^{|\cdot|_p}$ is a compact open subgroup of $\Que_p$.  
Also, by \cite[25.1]{hewittrI}, $\what{\Que_p}\cong\Que_p$.
If $\tau\in\tnq{\Que_p}$, then by Theorem \ref{theo:abdual}
$\hat{\tau}$ is a locally compact topology on $\Que_p$ which is finer than
$\tau_{\Que_p}$.  However, since $\Zee_p$ is an open subgroup,
it follows \cite[Theorem 2.1]{rickert}, that $\hat{\tau}$ is either 
$\tau_{\Que_p}$ or the discrete topology.  Thus $\tau=\tau_{\Que_p}$ or
$\tau=\tauap$, accordingly. \endpf

Returning to topologies on $\Que$, we let
\[
\tau_R=\til{\tau}_\Ree\vee\tauap,
\aand\tau_p=\til{\tau}_{\Que_p}\vee\tauap\text{ for
each }p\iin\Pee.
\]

\begin{padicandreal}\label{prop:padicandreal}
None of the topologies $\til{\tau}_\Ree$ or $\til{\tau}_{\Que_p}$ ($p\in\Pee$)
are non-quotient topologies on $\Que$.  However, $\tau_R$ and
$\tau_p$ ($p\in\Pee$) form distinct non-quotient topologies.
\end{padicandreal}

\proof First, let us note that $\what{\Ree}|_{\Que}\not=\what{\Que}$
by  \cite[25.26]{hewittrI}.  Hence it follows Theorem \ref{theo:abdual}
that $\til{\tau}_\Ree\not\in\tnq{\Que}$.  Moreover we have that
$\what{\Ree}|_\Que\cap\what{\Que_p}|_\Que=\{1\}$, the trivial
character.  Indeed, if $\sig\in\what{\Ree}\setdif\{1\}$ and $U$ is
a sufficiently small neighbourhood of the unit in $\Tee$, then
$\sig^{-1}(U)$ is a non-dense open subset of $\Ree$, and hence
cannot contain any of the basic neighbourhoods $r+p^kO_p$
from $\til{\tau}_{\Que_p}$.  It thus follows from Theorem \ref{theo:abdual}
that $\til{\tau}_{\Que_p}\not\in\tnq{\Que}$.

Now we wish to show that $\tau_R,\tau_p$ and $\tau_q$,
where $p,q$ are distinct primes, are distinct topologies.
Let $r_n=p^n/q^n\iin\Que$.  Then we have
\[
\til{\tau}_{\Que_q}\text{-}\lim_{n\to\infty}r_n=\infty\quad\text{ while }
\quad\til{\tau}_{\Que_p}\text{-}\lim_{n\to\infty}r_n=0.
\]
Hence $\{0\}\cup\{r_n\}_{n\in\En}$ is $\til{\tau}_{\Que_p}$-compact
and hence $\tau_p$-precompact, while it is $\til{\tau}_{\Que_q}$-unbounded
and hence not $\tau_q$-precompact.  Furthermore, if $p>q$ the same
argument shows $\tau_p\not=\tau_R$; and if $p<q$ we conduct
the same argument on $1/r_n$ to show that $\tau_q\not=\tau_R$.
\endpf

We can adapt the arguments above to show that the topologies
$\{\tau_R,\tau_p:p\in\Pee\}$ generate a free semilattice, 
$\fF=\langle\tauap,\tau_R,\tau_p:p\in\Pee\rangle$,
in $\tnq{\Que}$.  However, $\fF\not=\tnq{\Que}$.  
Consider, for example, the (open) subgroup $\Zee$.  If $\tauap^0$
denotes the topology induced by the almost periodic compactification
on $\Zee$, then we can extend this to a locally precompact topology
$\bar{\tau}$ on $\Que$, as in Lemma \ref{lem:exttop}. It can be shown
that $\bar{\tau}$ is not in $\fF$.

Let $\mathrm{G}(\Que)$ be any algebraic group which admits an
embedding $\Que\hookrightarrow\mathrm{G}(\Que)$.
An example
of such is $\mathrm{SL}_2(\Que)$ with the embedding
$r\mapsto\begin{pmatrix} 1 & r \\ 0 & 1\end{pmatrix}$.
Then if let $\tau_R$ and $\tau_p$ be the non-quotient
topologies whose respective quotients are those topologies
induced by the embeddings $\mathrm{G}(\Que)\hookrightarrow
\mathrm{G}(\Ree)$ and $\mathrm{G}(\Que)\hookrightarrow
\mathrm{G}(\Que_p)$.  Then these topologies must all
be distinct.  Thus we expect $\tnq{\mathrm{G}(\Que)}$
to be quite complicated.  On the basis of Propositions
\ref{prop:padic} and \ref{prop:reeandtee}, we conjecture that
$\tnq{\mathrm{G}(\Que_p)}\cong\tnq{\mathrm{G}(\Ree)}$.

\subsection{Free Groups}\label{ssec:freegrp}

Free groups will give us very complicated lattices of non-quotient
topologies.  To capture a glimpse of this, let $\mathrm{F}_\infty$
be the free group on generators $\{x_n\}_{n\in\En}$.  Then
if $G$ is any separable locally compact group and $S=\{s_n\}_{n\in\En}$ is
any set which generates a dense subgroup, then the map
\[
\eta_S:\mathrm{F}_\infty\to G\quad\text{ given by }\quad\eta_S(x_n)=s_n
\]
is a homomorphism.  We can thus generate different choices
for locally precompact topologies through different choices of
pairs $(G,S)$, chosen as above.  Thus, it seems to
be impossible to determine $\tnq{\mathrm{F}_\infty}$ in any conventional
terms.

Note that the above construction can be generalised to any discrete
free group $\mathrm{F}_X$ on any infinite set of generators.  Moreover, since
for each finite $n\geq 2$, the free group $\mathrm{F}_n$ contains an (open)
copy of $\mathrm{F}_\infty$, Lemma \ref{lem:exttop} shows that 
$\tnq{\mathrm{F}_n}$ is very complicated too.

A less tractable problem appears to be the determination of the non-quotient
topologies for the free abelian group $\Zee^{\oplus\infty}$.  
However pathologies still exist as
$\tnq{\Zee^{\oplus\infty}}$ is not a complete semilattice;  
the convenience of this property is noted in Proposition
\ref{prop:hdprincipal}.  Indeed, fix an irrational real number $\xi$, let 
\[
\eta_n:\Zee^{\oplus\infty}\cong\bigl(\Zee^2\bigr)^{\oplus\infty}\to
\Ree^n:(k_1,l_1,\dots,k_n,l_n,\dots)\mapsto(k_1+\xi l_1,\dots,k_n+\xi l_n)
\]
and let $\tau_n=\eta_n^{-1}(\tau_{\Ree^n})\vee\tauap$.   Then $\tau=
\bigvee_{n=1}^\infty\tau_n$ is a topology on $\Zee^{\oplus\infty}$
for which $\Zee^{\oplus\infty}_\tau$ admits $\Ree^{\oplus\infty}$
as a quotient group with compact kernel.  (Here
$\Ree^{\oplus\infty}=\bigcup_{n=1}^\infty \Ree^n\oplus 0^\infty$ 
has the inductive
limit topology.)  Hence $\tau$ is not locally precompact.

{
\bibliography{spinebib}
\bibliographystyle{plain}
}

\medskip
{\sc  Monica Ilie, Department of Mathematical Sciences, Lakehead University,
955 Oliver Road,  Thunder Bay, ON P7B\,5E1, Canada} 

E-mail: {\tt milie@lakeheadu.ca}

\medskip
{\sc Nico Spronk, Department of Pure Mathematics, University of Waterloo,
Waterloo, ON N2L\,3G1, Canada} 

E-mail: {\tt nspronk@uwaterloo.ca}

\end{document}